\RequirePackage{fix-cm}
\documentclass[twocolumn]{svjour3}          
\smartqed  
\usepackage{graphicx}
\usepackage{multicol}
\usepackage{graphicx}
\usepackage{xcolor}
\usepackage{amssymb, amsmath}
\usepackage{float}
\usepackage{enumitem}
\usepackage{mathtools, cuted}
\usepackage{verbatim}
\usepackage{tikz,xcolor}
\usepackage[hidelinks]{hyperref}
\usepackage{stfloats}
\usepackage{sidecap}
\usepackage{algorithm}
\usepackage[noend]{algpseudocode} 
\usepackage{array,ragged2e,multirow}
\usepackage{scrtime}
\usepackage[usestackEOL]{stackengine}
\usepackage{booktabs}
\usepackage[bb=dsserif]{mathalpha}

\newcommand{\window}[2]{\ensuremath{\mathit{window}_{#1, #2}}}
\newcommand{\machine}[2]{\ensuremath{\mathit{machine}_{#1, #2}}}
\newcommand{\sd}[1]{\ensuremath{\mathit{start\_day}_{#1}}}
\newcommand{\pinP}{p\in \mathcal{P}}
\newcommand{\pinPpref}{p\in \mathcal{P}^{\mathit{pref}}}
\newcommand{\dinD}{d \in \mathcal{D}}
\newcommand{\dinDw}{d \in \mathcal{D}_w}
\newcommand{\finF}{f \in \mathcal{F}}
\newcommand{\winW}{w \in \mathcal{W}}
\newcommand{\minM}{m \in \mathcal{M}}
\newcommand{\hinH}{h \in \mathcal{H}}
\newcommand{\iinK}{i \in \mathcal{K}_p}
\newcommand{\R}{\mathbb{R}}

\newcommand{\Dw}{\mathcal{D}_w}
\newcommand{\Bm}{\mathcal{B_M}}
\newcommand{\binBm}{b_\mathcal{M} \in \mathcal{B_M}}
\newcommand{\MprefP}{\mathcal{M}^{\mathit{pref}}_p}
\newcommand{\Ppref}{\mathcal{P}^{\mathit{pref}}}
\newcommand{\mgP}{\mathit{machine\_group}_p}

\newcolumntype{P}[1]{>{\RaggedRight\arraybackslash}p{#1}}

\newcommand*{\affaddr}[1]{#1} 
\newcommand*{\affmark}[1][*]{\textsuperscript{#1}}
\newcommand{\mypar}[1]{\vspace{-10pt} \paragraph{\textbf{#1.}}}

\definecolor{lime}{HTML}{A6CE39}
\DeclareRobustCommand{\orcidicon}{
	\begin{tikzpicture}
	\draw[lime, fill=lime] (0,0) 
	circle [radius=0.16] 
	node[white] {{\fontfamily{qag}\selectfont \tiny ID}};
	\draw[white, fill=white] (-0.0625,0.095) 
	circle [radius=0.007];
	\end{tikzpicture}
	\hspace{-3mm}
}

\foreach \x in {A, ..., Z}{%
	\expandafter\xdef\csname orcid\x\endcsname{\noexpand\href{https://orcid.org/\csname orcidauthor\x\endcsname}{\noexpand\orcidicon}}
}

\journalname{Journal}
\begin{document}

\title{Comparing Optimization Methods for Radiation Therapy Patient Scheduling using Different Objectives}

\titlerunning{Comparing Optimization Methods for Radiation Therapy Patient Scheduling...} 

\author{Sara Frimodig\affmark[1,2*]\orcidA{} \and
        Per Enqvist\affmark[1]\orcidB{} \and
        Mats Carlsson\affmark[3]\orcidC \and
        Carole Mercier\affmark[4]\orcidD{}
}

\authorrunning{S. Frimodig et al.}

\institute{   \affaddr{\affmark[1]Department of Mathematics, KTH Royal Institute of Technology, Stockholm, Sweden}  \\
              \affaddr{\affmark[2]RaySearch Laboratories, Stockholm, Sweden } \\
              \affaddr{\affmark[3]Department of Computer Science, RISE Research Institutes of Sweden }\\
              \affaddr{\affmark[4]Department of Radiation Oncology, Iridium Netwerk, Antwerp, Belgium}                 \\
            *Corresponding author, \email{sarhal@kth.com}
}

\date{Received: date / Accepted: date}

\maketitle

\begin{abstract}
Radiation therapy (RT) is a medical treatment to kill cancer cells or shrink tumors. To manually schedule patients for RT is a time-consuming and challenging task. By the use of optimization, patient schedules for RT can be created automatically. This paper presents a study of different optimization methods for modeling and solving the RT patient scheduling problem, which can be used as decision support when implementing an automatic scheduling algorithm in practice. We introduce an Integer Programming (IP) model, a column generation IP model (CG-IP), and a Constraint Programming model. Patients are scheduled on multiple machine types considering their priority for treatment, session duration and allowed machines. Expected future patient arrivals are included in the models as placeholder patients. Since different cancer centers can have different scheduling objectives, the models are compared using multiple objective functions, including minimizing waiting times, and maximizing the fulfillment of patients' preferences for treatment times. The test data is generated from historical data from Iridium Netwerk, Belgium's largest cancer center with 10 linear accelerators. The results demonstrate that the CG-IP model can solve all the different problem instances to a mean optimality gap of less than $1\%$ within one hour. The proposed methodology provides a tool for automated scheduling of RT treatments and can be generally applied to RT centers.

\end{abstract}

\section*{Statements and Declarations}

\textbf{Competing interests} The authors declare that they have no conflict of interest.

\noindent \textbf{Availability of data and material} The datasets generated during the current study are available in the Open Science Framework repository, and are accessible through \href{https://osf.io/45qw2/?view_only=4a0a67e21cb542df8f9a0f74241de825}{\textit{this link}.} 



\begin{acknowledgements}
The authors want to extend thanks to Kjell Eriksson, Chief Science Officer at RaySearch Laboratories, for valuable discussions, and Geert De Kerf, Medical Physics Expert at Iridium Netwerk, for help collecting the data used in the experiments and providing necessary clinical information. 
\end{acknowledgements}

\clearpage
\section{Introduction}
Radiation therapy (RT) is a cancer treatment that uses radiation to kill malignant tumor cells. Together with chemotherapy and surgery, it is one of the most commonly used cancer therapies worldwide. Based on demographic changes such as an aging population, cancer incidents are increasing, and a 16\% expected increase in the number of RT treatment courses in Europe has been estimated from 2012 to 2025 \cite{Borras2016}. 

A long waiting time between referral and treatment start has negative effects on the outcome of the treatment. Reasons for this include tumor growth, psychological distress of the patient, and prolonged symptoms when the waiting times are long \cite{CHEN2008,Fortin2002,Gomez2015,ORourke2000,VanHarten2015}. 
Therefore, many cancer institutes around the world have adopted waiting time targets that state the maximum allowed waiting time before treatment starts. 

The RT treatment intent can be divided into either curative or palliative, where the first intends to cure the patient, and the latter mainly aims to provide symptom relief. 
Furthermore, cancer patients are often divided into different urgency levels depending on the site of the cancer, treatment intent, and the size and progress of the tumor. The prioritization of patients for treatment can be done in different ways in different countries or hospitals  \cite{Lim2005,Scoccianti2012,Ebert2013,Thomsen2009}.

There are several types of RT, where external photon beam RT is by far the most common, and the one covered in this paper. Photon beam RT is delivered on machines called linear accelerators (\textit{linacs}). Because the DNA of healthy cells is repaired to a higher degree than that of malignant cells, the radiotherapy treatment is usually divided into multiple sessions, called \textit{fractions}. The fractions are scheduled daily with breaks on the weekends for a period of up to eight weeks. The duration of the fractions varies between patients due to treatment technique and treatment complexity. For a particular treatment, the delivery time of all fractions is the same. However, the first fraction is usually scheduled for a longer time since it includes setup times, extra time for patient education and reassurance, and additional quality checks before treatment start \cite{DeVitaJr2018}.

In the RT patient scheduling problem, the aim is to schedule RT treatments for a set of patients, given a set of linacs, for a certain planning horizon. The problem is complicated since the patients are of different priorities, and their treatments vary in the number of fractions, fraction durations and set of compatible linacs. Furthermore, the RT process includes many uncertainties, such as the random arrival of new patients. 

This paper considers the RT patient scheduling problem arising at Iridium Netwerk, an RT center located in Antwerp, Belgium. In 2020, they operated 10 linacs, delivering 5500 RT treatments to approximately 4000 patients. The scheduling at Iridium is today done manually. Designing more efficient schedules would be of great significance as it could potentially improve patient outcomes by shortening waiting times. For this reason, this paper makes the following contributions:
\begin{itemize}
    \item The \emph{main contribution} is a comprehensive comparison and performance study of three exact optimization approaches to the RT scheduling problem. By evaluating the suitability of different methods for the problem, this serves as a foundation for further research in the area, and more importantly, as a decision support when implementing an automatic scheduling algorithm in practice. 
    \item The main \emph{technical novelty} lies in the original models developed: an integer linear programming (IP) model, a column generation IP (CG-IP) model, and a constraint programming (CP) model, as well as a method combining the CP and IP models. 
    To the best of our knowledge, these models are the first to simultaneously assign all fractions of the patients to both linacs and specific time windows, while including all the medical and technical constraints necessary for the scheduling to work in practice. 
    Furthermore, it is the first time column generation has been used for the RT patient scheduling problem. The problem instances solved are larger than in previous studies in terms of number of linacs, number of patients and length of the planning horizon. 
    \item Different cancer centers may have different intentions when creating the RT schedules. In order to study the suitability of the above-mentioned optimization methods for various cancer centers, each model is solved using multiple different objective functions to evaluate if some particular optimization method is better for a certain objective. The objectives include minimizing waiting times and maximizing the satisfaction of time window preferences among the patients. 
\end{itemize}
The paper is organized as follows. Section~\ref{B-Sec:rel_work} presents the related work. Section~\ref{B-Sec: Problem Formulation} describes the problem. Section~\ref{B-Sec: Models} presents the models. The setup for the computational experiments is explained in Section \ref{B-Sec: Experiments}, followed by the results in Section \ref{B-Sec: Results}. Section \ref{B-Sec:Discussion} contains the discussion, and Section \ref{B-Sec: Conclusions} presents the conclusions.

\section{Related Work}
\label{B-Sec:rel_work}
A literature review on the use of operations research for resource planning in RT was published in 2016 by Vieira et al. \cite{Vieira2016}. The authors found $12$~papers addressing the problem of scheduling RT patients on linacs. The first use of integer programming (IP) for optimization of RT appointments was published in 2008 by Conforti et al. \cite{Conforti2008}, where a block scheduling model is presented. The day is divided into blocks of equal duration and each treatment is assigned to one block. The same authors later developed a non-block scheduling model, which allows for different treatment durations \cite{Conforti2010}. Another IP model for non-block scheduling is presented by Jacquemin et al. \cite{Jacquemin2011}, where the notion of treatment patterns is introduced to allow non-consecutive treatment days. A limitation of these papers is that they do not consider all the constraints present in RT scheduling, such as multiple machine types and partial availability in the schedule. 

Sauré et al. \cite{Saure2012} present a method for advance RT patient scheduling using a discounted infinite-horizon Markov decision process, and show that their proposed policy can increase the percentage of treatments initiated within 10 days from 73\% to 96\%. In \cite{Gocgun2018}, Gocgun  extends the same problem setup to also include patient cancellations. The setting used in these papers is a simplified model of a cancer center, equipped with three identical machines and 18 treatment types. The resulting policies assign a start day to each patient, with no sequencing of patients throughout the day. 

In order to schedule RT appointments one by one in an online fashion, Legrain et al. \cite{Legrain2015} propose a hybrid method combining stochastic and online optimization using a block-scheduling strategy. The results show that their method works well on small real instances, with two linacs and less than 3.5 requests per day. Aringhieri et al. \cite{Aringhieri2020} also present methods for online RT scheduling, and develop three online optimization algorithms for a block-scheduling formulation and one machine.

Li et al. \cite{Li2019} model the RT patient scheduling as a queueing system with multiple queues. A new class of scheduling policies is proposed, where the parameters are selected through simulation-based optimization heuristics. All treatments are assumed to have the same length, and all machines are assumed to be identical. 

The type of RT that uses high-energy particles (protons or helium ions) is referred to as particle therapy (PT), in which a single particle beam is shared between multiple treatment rooms. This gives very different technical and medical constraints than in conventional photon beam RT. Two papers that present methods for optimizing the sequencing of patients throughout the day in PT are Vogl et al. \cite{Vogl2019} and Maschler et al. \cite{Maschler2020}, both aiming to schedule treatments close to a pre-defined target time and both using different heuristic methods. Braune et al. \cite{Braune2021} present a model for planning appointment times in PT under uncertain activity durations, and solve the resulting stochastic optimization model using a combination of a Genetic Algorithm and Monte Carlo simulation. 

The first paper to include the sequencing of patient throughout the day in photon beam RT was published in 2020 by Vieira et al. \cite{Vieira2020}. The authors create weekly schedules with the objective to maximize the fulfillment of the patients' time window preferences using a mixed-integer programming (MIP) model together with a pre-processing heuristic to divide the problem into subproblems for clusters of machines. In a second paper they test their method in two Dutch clinics \cite{Vieira2021}, with results showing that the weekly schedule was improved in both centers. 
However, the problem studied in these papers is different from the one in this paper, as they create weekly schedules for a time horizon of five days, whereas we aim to assign all fractions to machines and therefore have a significantly longer time horizon (typically around 80 days). 
Using a data-driven approach,  Moradi et al. \cite{Moradi2022} study the patient sequencing problem in a simplified clinical setup, where all treatment durations are equal and all machines are identical and independent. To improve the weekly schedules, the authors utilize patient information in a MIP model to determine the optimal sequence of patients for a list of patients that have been previously assigned to a treatment day. The results show that it is favorable to schedule reliable patients early on to reduce idle time on machines caused by delayed patients or no-shows.

Constraint Programming (CP) is a technique for solving combinatorial problems with origins in the computer science and artificial intelligence communities. CP solvers classically explore the solution space using tree search-based heuristics. For an overview, see~\cite{Rossi2006}. CP has been used in RT treatment planning~\cite{Baatar2011}, in chemotherapy patient scheduling~\cite{Hahn-Goldberg2014} and in operating room scheduling~\cite{Hashemi2016}. Overall, scheduling is a field where CP has shown to be effective, see for example~\cite{Bartak2010}. Frimodig et al. \cite{Frimodig2019} present and compare two CP models and one IP model for the RT scheduling problem. The CP models are shown to be efficient at finding feasible solutions, but are generally slower than the IP model at proving optimality. 
A limitation for these models is that they only consider one machine type. 

Pham et al. \cite{Pham2021} propose a two-stage approach for the RT scheduling problem. In a first phase, an IP model is used to assign patients to linacs and days, and in the second phase the patients' specific appointment times are decided using either a MIP or a CP model. The test data is generated based on data from CHUM, a cancer center in Canada. The test instances have seven linacs and a time horizon of 60 days. The results in the second phase show that CP finds good solutions faster, while MIP is better at closing optimality gaps, which agrees with the results in \cite{Frimodig2019}. 
Some simplifications in their models are that all patients can be treated on all machines, that all machine switches are allowed, and that the first fraction has the same duration as the rest. These assumptions make the models less general than the ones presented in this paper, and not suited for the scheduling problem at Iridium Netwerk. 

Column generation (CG) is a method that is often successful in solving certain classes of large scale integer programs. The method alternates between a restricted master problem and a column generation subproblem. CG has been applied in various areas within the medical treatment field, such as for surgeon and surgery scheduling \cite{Wang2018}, for patient admission \cite{Range2014}, and for nurse scheduling \cite{Bard2005}. In RT, it has been used for brachytherapy scheduling using deteriorating treatment times \cite{Shao2021,Shao2023}. In brachytherapy, the radiation is produced from a radioactive source placed within the patient, and the problem differs significantly from patient scheduling in conventional RT.

\section{Problem Formulation}
\label{B-Sec: Problem Formulation}

The RT scheduling problem consists of assigning each fraction for each patient to a day, a time window, and a machine. This section presents the real-world constraints and objectives present at Iridium Netwerk. 

\begin{figure*}[b]
\begin{center}
    \includegraphics[width=0.85\linewidth]{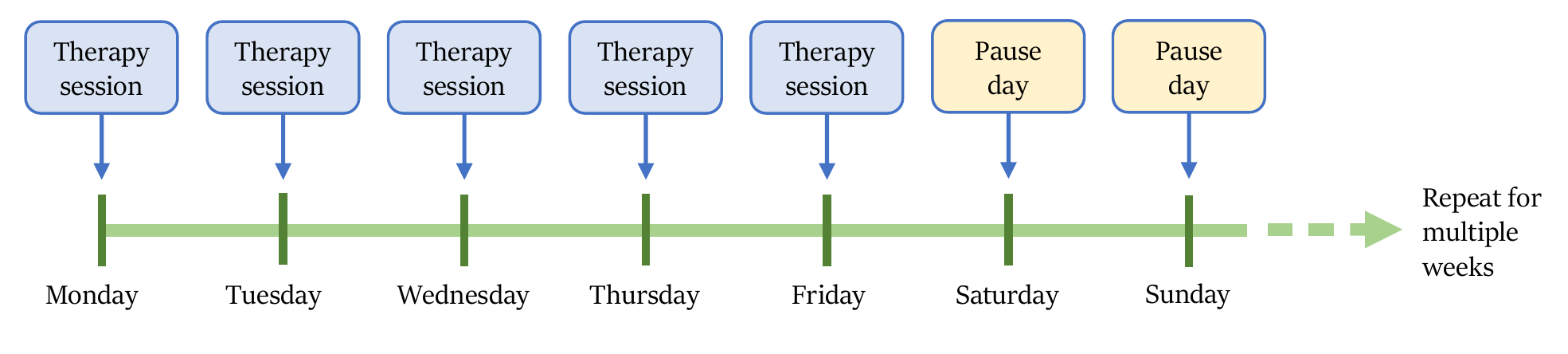}
    \caption{A typical fractionation scheme of an RT patient}
 \end{center}
\label{B-fig:treatment_course}
\end{figure*}

\mypar{Patients}
A \emph{priority} is assigned to each patient based on urgency and treatment intent. The prioritization can be done differently in different countries or hospitals \cite{Lim2005,Scoccianti2012,Ebert2013,Thomsen2009}, but at Iridium Netwerk it is done by a physician. In 2020 at Iridium, there were three priority groups, and approximately $42\%$ patients were priority~A, $18\%$ priority~B, and $40\%$ priority~C. Furthermore, each patient is assigned to a \emph{treatment protocol}, which states the \emph{fractionation scheme} (that is, how many days the patient is to be treated and with which frequency), and the \emph{duration} of the first and subsequent fractions. An example of a fractionation scheme is shown in Figure \ref{B-fig:treatment_course}. Different protocols have different allowed \emph{start days}: palliative patients can start any weekday, whereas curative patients cannot start on Fridays. Both the number of priority groups and the allowed start days are easily generalizable in the models.  Examples of some different treatment protocols can be seen in Table \ref{B-tab:protocols}.

Urgent patients must start treatment soon after arrival. Since the patients are of different priority groups, and the fractionation schemes typically span multiple weeks, this must be considered when creating the schedules. In practice, most clinics handle this by reserving empty time slots on each machine for urgent patients. In this paper, the expected value of the future patient arrivals for the coming weeks are included in the models as placeholder patients (i.e., \textit{dummy} patients) to predict the expected utilization of resources. The models are deterministic, and the placeholder patients are handled as regular patients in the patient list.

\begin{table*}
\caption{Examples of some different treatment protocols at Iridium Netwerk}
\label{B-tab:protocols}
\begin{tabular}{P{2.82cm} P{1.08cm} P{1.73cm} P{2.06cm} P{1.2cm}  P{3.05cm} P{2.5cm}}
\hline\noalign{\smallskip}
Protocol & Priority & Duration first fraction (minutes) & Duration other fractions (minutes) & Average number fractions& Preferred machines &Allowed machines\\
\hline\noalign{\smallskip}
Bladder (VMAT) & B & 24 & 12 & 20 & M2, M3, M5, M6, M7, M9 & M1, M4, M8\\
Brain STX 3x & A & 40 & 40 & 3 & M10 \\
Breast bilateral & C & 48 & 24 & 13 & M1, M3, M4, M5, M6, M8 & M2, M7 \\
Head-Neck (VMAT) & A & 24 & 12 & 26 & M2, M3, M5, M6, M9 & M1, M4 \\
Prostate SBRT & C & 24 & 12 & 5 & M10 \\
Rectum 25x & B & 24 & 12 & 25  & M2, M3, M5, M6, M8 & M1, M4, M7 \\
\hline\noalign{\smallskip}
\end{tabular}
\end{table*}

\mypar{Time}
RT clinics have different routines for scheduling patients on linacs. Some gather patients into a \emph{batch} and schedule them once or several times per day, while others schedule each patient at admission. In \cite{Pham2021}, different scheduling strategies are evaluated using a simulation. Preliminary results show that daily batch scheduling reduces patients’ waiting time and overdue time. This paper focuses on batch scheduling and assumes that the scheduling is done at the end of each day taking patients from previous days into account.

There are two different scheduling systems used at RT centers: \emph{block} or \emph{non-block} scheduling systems. The block system, where the day is divided into blocks of equal duration and each patient is assigned to a block, is more widely used in clinics, but has severe drawbacks since there is no way to control the variability of treatment time. This can generate costs related to machine under utilization, staff overtime, and patient waiting time. This paper uses a \emph{non-block} scheduling strategy, where a day is  instead divided into \emph{time windows}. A time window is typically $1.5-4$~hours while a treatment duration is $10-60$~minutes, which makes this is different from a block scheduling system where each treatment is assumed to have the same duration as one or multiple blocks. Patients are assigned to windows instead of specific start times, as this leads to simpler and more efficient models while maintaining an adequate level of detail from a clinical perspective. The specific start time for each treatment within the time window is given in a post-processing step.

During the first treatment fraction, extra time is needed for both instructing the patient and for linac setup \cite{Turner2002}. Therefore, auxiliary time must be assigned to each new patient, which is done by assigning the first fraction to a longer time duration (see Table \ref{B-tab:protocols}). Furthermore, at Iridium no patients are treated during weekends. This fact is used to simplify the models; the time horizon is adjusted to only contain weekdays ($\mathcal{D}_w$). In general, at most cancer centers in the world only patients undergoing an oncologic emergency are treated during weekends, and in that case, the care is not planned more than a day in advance \cite{Mitera2009,Yeo2012}.

\mypar{Machines}
The radiation is delivered on linacs. This paper assumes that there are multiple machine types, which is the case in almost all clinics. At Iridium there are ten linacs distributed over four different hospitals. The treatment protocol for each patient states one or more machines that can deliver the protocol, however, some machines are preferred over the others. An example can be seen in Table \ref{B-tab:protocols}.

Some linacs are so called \textit{beam-matched}, meaning that a patient can switch between these linacs between fractions. Two linacs are considered completely beam-matched if they are the same machine type at the same hospital, and partially beam-matched if they are the same machine type, but at different hospitals. Switching between completely beam-matched machines can be done at no cost, whereas there is a cost for switching to a machine that is only a partially matched. The beam-matched machines are presented in Table \ref{B-tab:beam-matched}. To generalize the models, the cost for  machine switches between partially beam-matched machines is not active in all objective functions investigated.
\begin{table}
\caption{Beam-matches machines at Iridium Netwerk}
\label{B-tab:beam-matched}
\begin{tabular}{P{1.15cm} P{2.9cm} P{3.1cm}}
\hline\noalign{\smallskip}
Machine & Completely matched & Partially matched\\
\hline\noalign{\smallskip}
M1 & - & M4, M8 \\
M2 & - & M3, M5, M6, M7, M9 \\
M3 & M9 & M2, M5, M6, M7 \\
M4 & - & M1, M8 \\
M5 & M6 & M2, M3, M7, M9 \\
M6 & M5 & M2, M3, M7, M9 \\
M7 & - & M2, M3, M5, M6, M9 \\
M8 & - & M1, M4 \\
M9 & M3 & M2, M5, M6, M7 \\
M10 & - & - \\
\hline\noalign{\smallskip}
\end{tabular}
\end{table}

\mypar{Objectives}
When creating the RT schedules, different cancer centers can have different scheduling objectives. In order to compare the suitability of different optimization methods to solve the RT scheduling problem, the models are tested for multiple objective functions to evaluate their performance, and also to evaluate if some particular optimization method is better for a certain scheduling objective. The following objectives will be tested in different combinations:
\begin{enumerate}[label=(\roman*)]
\item \label{B-enum:i} Minimize a weighted sum of the waiting times
\item \label{B-enum:ii} Minimize a weighted sum of the violations of the target dates
\item \label{B-enum:iii} Minimize the number of time window switches
\item \label{B-enum:iv} Minimize violations of time window preferences 
\item \label{B-enum:5} Minimize the number of fractions scheduled on non-preferred machines
\item \label{B-enum:6} Minimize the number of switches between machines that are not completely beam-matched
\end{enumerate}
The weights in the weighted sums in \ref{B-enum:i} and \ref{B-enum:ii} should reflect the severeness of delaying treatment start for the different priority groups. In objective \ref{B-enum:ii}, the waiting time targets are assumed to be $2$~days for priority~A, $14$~days for priority~B, and  $28$~days for priority~C patients, but this is easily generalizable. The waiting time targets differ between countries, and advanced methods to determine the waiting time targets have recently been studied \cite{Babashov2023}. The aim for objective \ref{B-enum:iii} is to schedule patients at approximately the same time each day, since this is something the booking administrators usually try to do. Literature shows that patients have different preferences regarding the time of their appointments \cite{Olivotto2015}, which is what should be captured in objective \ref{B-enum:iv}. As many fractions as possible should be scheduled on the machines preferred for the particular treatment, hence objective \ref{B-enum:5} states that the number of fractions scheduled on a non-preferred machine should be minimized. Finally in objective \ref{B-enum:6}, the number of switches between machines that are not completely beam-matched should be minimized. In Section \ref{B-sec:objective_functions}, the  objectives will be combined into different objective functions. For example, a combination of \ref{B-enum:i}-\ref{B-enum:6}, with \ref{B-enum:i} being most important, is most similar to what is used at Iridium Netwerk.

\begin{table*}[t]
\caption{Notations for the models}
\label{B-tab:inputs}
\begin{tabular}{P{0.18\linewidth} P{0.77\linewidth}}
\hline\noalign{\smallskip}
    Parameter                                  & Description \\
    \hline\noalign{\smallskip}
	$\mathcal{P} = \{1, \dots, P\}$            & Set of all patients, $P \in \mathbb{N} $   \\
	$\Dw = \{1, \dots, D_w \}$                      & Set of weekdays. $D_w \in \mathbb{N}$ is the number of weekdays in the planning horizon \\
	$\mathcal{W} = \{1, \dots, W\}$            & Set of time windows in a day, $W \in \mathbb{N}$  \\
	$L_w \in \mathbb{N}$                       & The window length for window $\winW$ in number of minutes\\
	$\mathcal{M} = \{1, \dots, M\}$            & Set of machines, $M \in \mathbb{N} $  \\
	$\mathcal{M}_p \subseteq \mathcal{M}$      & Set of machines allowed for patient $p$  \\
	$\MprefP \subseteq \mathcal{M}_p$          & Set of machines preferred for patient $p$  \\
	$\mathcal{C_M}$                            & List of sets of completely beam-matched machines\\
	$\mathcal{P_M}$                            & List of sets of partially beam-matched machines\\
	$\Bm = \mathcal{C_M} \cup \mathcal{P_M}$   & List of sets of all beam-matched machines\\
	$\mathcal{H}=\{1,\dots,H\} $               & Set of treatment protocols, $H=72$ \\
    $\mathit{dur}_{p0} \in \mathbb{N}$         & Duration of first fraction for patient $p$ (minutes) \\
	$\mathit{dur}_p \in \mathbb{N}$            & Duration of fractions other than first for patient $p$ (minutes) \\
	$\mathcal{F}_p \in \{1,\dots, F_{p} \}$     & Set of all fractions for patient $p$. $F_p \in \mathbb{N}$ is the number of fractions 
	\\
    $S \in \R^M \times \R^D \times \R^{W}$   & The number of occupied timeslots in each window, machine and day, i.e. $S_{m,d,w} \in \{0,...,L_w\}$ \\
	$\mathcal{A}_p \in \Dw$                  & The set of allowed start days for the protocol of patient $p$ \\
	$c_p \in \{10,3,1\}$                       & Weights for patient $p$ in priority group A, B or C \\
	$d_{L,p} \in \mathcal{D}_w$                & The day limit for latest allowed treatment start for patient $p$, adjusted for days already waited \\
	$d_{\min,p} \in \mathcal{D}_w$             & The earliest day for patient $p$ to be scheduled \\
	$\Ppref \subset \mathcal{P}$   & The set of patients that have a time window preference \\
	$w^{\mathit{pref}}_{p} \in \mathcal{W}$   & The window preference of patient $\pinPpref$ \\
\noalign{\smallskip}\hline
\end{tabular}
\end{table*}

\begin{table}
\caption{Variables in the IP model}
\label{B-tab:variables_IP}
\begin{tabular}{P{0.29\linewidth} P{0.61\linewidth}}
\hline\noalign{\smallskip}
	$q_{p,m,d,f} \in \{0,1\}$                 & $1$ if patient $\pinP$ has their $f$th fraction ($\finF_p$) on weekday $\dinDw$ on machine $\minM$, $0$ otherwise\\
		$x_{p,m,d,w} \in \{0,1\}$             & $1$ if patient $\pinP$ is scheduled in window $\winW$ on machine $\minM$ on weekday $\dinDw$, $0$ otherwise \\
		$t_{p,m,d,w} \in \{0,1\}$             & $1$ if patient $\pinP$ starts treatment in window $\winW$ on machine $\minM$ on weekday $\dinDw$, $0$ otherwise \\
	$y_{p,d,w} \in \{0,1\}$                   & $1$ if patient $\pinP$ is scheduled in window $\winW$ on day $\dinDw$ \\
	$z_{p,d} \in \{0,1\}$                       & $1$ if patient $\pinP$ has switched windows from day $d$ to $d+1$, $0$ otherwise\\
	$u_{p,d} \in \{0,\dots,W-1\}$           & The violation of the time window preference for patient $\pinP$ on day $\dinD$    \\
	$v_{p,f} \in \{0,1\}$                   & $1$ if patient $\pinP$ is scheduled on non-preferred machine on fraction $\finF_p$, $0$ otherwise   \\
	$s_{p,f} \in \{0,1\}$                   & $1$ if patient $\pinP$ switches to a partially beam-matched machine from fraction $f$ to $f+1$, $0$ otherwise   \\
\noalign{\smallskip}\hline
\end{tabular}
\end{table}

\section{Models}
\label{B-Sec: Models}
Multiple models are developed: an IP model, a CG-IP model, and a CP model, as well as a method combining CP and IP. They are designed to capture the same real-world constraints and objectives. Their inputs are presented in Table \ref{B-tab:inputs}. As stated in Section \ref{B-Sec: Problem Formulation}, no patients are treated during weekends. Therefore, the time horizon is adjusted so that $\mathcal{D}_w$ only contains weekdays.

\subsection{Integer Programming Model}
\label{B-Sec:IP pure}
The variables in the IP model are presented in Table \ref{B-tab:variables_IP} and the formulation is stated in \eqref{B-Eq:IP objective}-\eqref{B-Eq:MIP_integers}. 

\begin{figure*}[t]
\begin{alignat}{3}
	 & \text{minimize} \quad        && 
	 \mathmakebox[0pt][l]{\begin{multlined}1+ \sum_{\pinP} (\alpha_1 f_{1,p}+ \alpha_2 f_{2,p} + \alpha_3 f_{3,p} +\alpha_4 f_{4,p} + \alpha_5 f_{5,p} + \alpha_6 f_{6,p})
	 \end{multlined}} \label{B-Eq:IP objective} &&\\
	 & \begin{multlined}\text{subject to} \\\end{multlined}  \quad &  &  \begin{multlined}\sum_{m \in b_\mathcal{M}} q_{p,m,d,f} = \sum_{m \in b_\mathcal{M}} q_{p,m,d+1,f+1}, \end{multlined} \enskip   &  & \begin{multlined}\forall \pinP, \binBm,\\ d=\{1,\dots,D_w-1\}, f = \{1,\dots,F_p-1\}\end{multlined} \label{B-Eq:q(p,m,k,f)=q(p,m,k+1,f+1)} \\
	 &                   &  & \sum_{\minM} \sum_{\dinDw} q_{p,m,d,f} = 1, \enskip                                                                                                                                                            &  & \forall \pinP, \finF_{p}                                       \label{B-Eq:sum(q) over d,m =1}           \\ 
	 &                   &  & q_{p,m,d,1} = \sum_{\winW} t_{p,m,d,w} , \enskip                                                                                                                                                            &  & \forall \pinP, \minM, \dinDw                                   \label{B-Eq:q_1 = t}           \\ 
	 &                   &  & t_{p,m,d,w} \leq x_{p,m,d,w} , \enskip                                                                                                                                                            &  & \forall \pinP, \minM, \dinDw  , \winW                                 \label{B-Eq:t <= x}           \\ 
	 &                   &  &\begin{multlined} q_{p,m,d,1}  = 0, \\ \end{multlined} \enskip                                                                                                                                                                                     &  & \begin{multlined}\forall  \pinP, \minM, \dinDw \text{ if } d > d_{\max,p}\\ \text{ or } d<d_{\min,p} \text{ or } d \notin \mathcal{A}_p\end{multlined} \label{B-Eq:q = 0 daylim}         \\
	 &                   &  &\begin{multlined} q_{p,m,d,f}  = 0, \\ \end{multlined} \enskip                                                                                                                                                                                     &  & \begin{multlined}\forall \pinP, \minM \text{ if } m \notin \mathcal{M}_p, \dinDw,\\ \finF_p \text{ if } f \notin \mathcal{F}_{p,d}  \end{multlined}        \label{B-Eq:q = 0 M_p}          \\
	 &                   &  & \sum_{\winW} x_{p,m,d,w} = \sum_{\finF_p} q_{p,m,d,f}, \enskip                                                                                                                                                  &  & \forall \pinP, \minM, \dinDw                                      \label{B-Eq:sum(x) = sum(q)}         \\
	 &                   &  & \begin{multlined}\sum_{\pinP} \Big((x_{p,m,d,w} - t_{p,m,d,w}) \mathit{dur}_p +\\[-1ex] \quad  t_{p,m,d,w} \mathit{dur}_{p0}\Big) + S_{m,d,w} \leq L_w,\end{multlined}  \enskip                                                                                                              &  & \forall \minM, \dinDw, \winW                                      \label{B-Eq:IP x*dur+S<=L_w}          \\
    &                    &  &	  \sum_{\dinDw} d \sum_{\minM}q_{p,m,d,1} \leq  \sum_{\dinDw}d \sum_{\minM} q_{p+1,m,d,1}, \enskip & & \forall h \in \mathcal{H}, p \in \mathcal{P}_h \text{ where } d_{L,p} \leq d_{L,p+1} \label{B-Eq:MIP1_dominance}\\
	 &                   &  & y_{p,d,w} \geq \sum_{\minM} x_{p,m,d,w}, \enskip                                                                                                                                                              &  & \forall \pinP, \dinDw, \winW                                      \label{B-Eq:MIP_y_1}               \\
	 &                   &  & \sum_{\winW} y_{p,d,w} = 1, \enskip                                                                                                                                                              &  & \forall \pinP, \dinDw                                              \label{B-Eq:MIP_y_2}                \\
	 &                    &  & z_{p,d} \geq y_{p,d,w}-y_{p,d+1,w}, \enskip                                                                                                                                                                &  & \forall \pinP, d=\{1,\dots,D_w-1\}, \winW                                     \label{B-Eq:MIP_z_1}\\
	 &                   &  & z_{p,d} \geq y_{p,d+1,w}-y_{p,d,w}, \enskip                                                                                                                                                           &  & \forall \pinP, d=\{1,\dots,D_w-1\}, \winW                                       \label{B-Eq:MIP_z_2} \\
	 &                   &  & u_{p,d} \geq \sum_{\minM} \sum_{\winW} x_{p,m,d,w}|w-w^{\mathit{pref}}_{p}|, \enskip                                                                                                                                                           &  & \forall \pinPpref, \dinD,                                       \label{B-Eq:MIP_u1} \\
	 &                   &  & u_{p,d} = 0, \enskip                                                                                                                                                           &  & \forall p \notin \pinPpref,  \dinD,                                      \label{B-Eq:MIP_u2} \\
 	 &                    &  & s_{p,f} \geq \sum_{\dinD} \sum_{m \in c_\mathcal{M}} (q_{p,m, d, f} - q_{p,m, d, f+1}), \enskip && \forall \pinP, c_\mathcal{M} \in \mathcal{C_M}, f = \{1,\dots,F_p-1\} \label{B-Eq:IP_machineswitch}\\
	&                    &  & v_{p,f} = \sum_{\minM} \sum_{\dinD} q_{p,m,d,f} \mathbb{1}_{(m \notin \mathcal{M}^{\mathit{pref}}_p)}, \enskip && \forall \pinP, \finF_p \label{B-Eq:IP_machinePref}\\
	 &                   &  &  \mathmakebox[0pt][l]{x\in\{0,1\}, q \in \{0,1\}, y \in \{0,1\}, z\in \{0,1\}, u \text{ integer}, s \in \{0,1\}, v \in \{0,1\}}&  & \label{B-Eq:MIP_integers}
\end{alignat}
\end{figure*}

\mypar{Constraints}
Constraint \eqref{B-Eq:q(p,m,k,f)=q(p,m,k+1,f+1)} is formulated to ensure that all fractions are scheduled after each other, and that they are all scheduled on beam-matched machines. Constraint \eqref{B-Eq:sum(q) over d,m =1} forces the $f$th fraction to be scheduled exactly one time for each patient. Constraint \eqref{B-Eq:q_1 = t} states that the first fraction for patient $p$ is scheduled on machine $m$ on day $d$, in any window, whereas Constraint \eqref{B-Eq:t <= x} also gives the correct window $w$ for the first fraction. 

The earliest day to start treatment is $d_{\min,p}$, and the last day to start treatment is $d_{\max,p} = D_w-F_p+1$.  A treatment can only start on an allowed start day given by $\mathcal{A}_p$. Furthermore, patient $p$ can only be scheduled on a machine that is allowed for the patient protocol given by $\mathcal{M}_p$. Finally, the set of allowed day-fraction pairs for patient $p$ is denoted $\mathcal{F}_{p,d} = \{\dinDw, \finF_p : f<d \text{ and } d-f < D_w-F_p \}$ (e.g., cannot schedule fraction 2 on day 1, or fraction 3 on day 50 if $F_p>3$ and the planning horizon $D_w=50$). In total, this is captured in Constraints \eqref{B-Eq:q = 0 daylim} and \eqref{B-Eq:q = 0 M_p}.

Constraint \eqref{B-Eq:sum(x) = sum(q)} states that each patient is scheduled in exactly one time window for each fraction. Constraint \eqref{B-Eq:IP x*dur+S<=L_w} ensures that all treatments fit within each time window.  The first term sums the session duration of all patients scheduled in window $w$ on machine $m$ on day $d$, except if it is the first fraction (since it will then evaluate to zero). The second term sums the durations of first fractions for patients starting their treatment in window $w$ on machine $m$ on day $d$. 
The sum of all scheduled patients' durations plus the already occupied time slots $S_{m,d,w}$ in that window should be less than or equal the window length $L_w$. 

For two patients with the same treatment protocol, the one with shorter day limit should start treatment first. Constraint \eqref{B-Eq:MIP1_dominance} enforces this 
by multiplying the variable $q_{p,m,d,1}$  with the day to get the start day, and force the ordering of the start days. Note the abuse of notation, where $p+1$ denotes the next entry in $\mathcal{P}_h$.

\mypar{Objective function} The objective functions \ref{B-enum:i}-\ref{B-enum:6} presented in Section \ref{B-Sec: Problem Formulation} are formulated. An offset set to $1$ is included to enable computation of the relative gap also when the optimal value is zero. The different objectives are combined with weights $\alpha_1,\dots, \alpha_6$ in \eqref{B-Eq:IP objective}.

Objective \ref{B-enum:i} is to minimize a weighted sum of the waiting times, which is formulated in \eqref{B-Eq:f1}. The number of waiting days after $d_{\min,p}$, the earliest day to be scheduled for patient $p$, are linearly penalized with weight $c_p$ corresponding to the priority group of patient $p$.
\begin{align}
f_{1,p} =  c_p \sum_{\minM} \sum_{d=d_{\min,p}}^{D_w} q_{p,m,d,1}(d-d_{\min,p})
\label{B-Eq:f1}
\end{align}
Objective \ref{B-enum:ii} is to minimize a weighted sum of the violations of the waiting time targets, formulated in objective \eqref{B-Eq:f2}. The days past the waiting time target $d_{L,p}$ are linearly penalized with weight $c_p$.
\begin{align}
f_{2,p} =  c_p \sum_{\minM} \sum_{d=d_{L,p}}^{D_w} q_{p,m,d,1}(d-d_{L,p})
\label{B-Eq:f2}
\end{align}
Objective \ref{B-enum:iii} is to minimize the number of time window switches for each patient. Therefore, the variable $y_{p,d,w}$ is defined according to Equations \eqref{B-Eq:MIP_y_1} and \eqref{B-Eq:MIP_y_2}, stating that if patient $p$ is scheduled on day $d$, $y_{p,d,w}=1$ only for the window where $p$ is scheduled. Every time window switch between two days is computed by  $|y_{p,d,w}-y_{p,d+1,w}| \quad \forall \pinP, \dinDw, \winW$. To avoid having the absolute value in the objective function, the variable $z_{p,d}$ is instead defined according to Constraints \eqref{B-Eq:MIP_z_1} and \eqref{B-Eq:MIP_z_2} and used in the objective function \eqref{B-Eq:f3}. 
\begin{align}
f_{3,p} = \sum_{\dinDw} z_{p,d}
\label{B-Eq:f3}
\end{align}
To form the objective corresponding to \ref{B-enum:iv}, the variable $u_{p,d}$ is defined by Constraints \eqref{B-Eq:MIP_u1} and \eqref{B-Eq:MIP_u2}. The time preference violation is zero if the patient does not have a preference, and is otherwise measured by the deviation from the preference on each day the patient is scheduled. Summing all violations gives \eqref{B-Eq:f4}.
\begin{align}
f_{4,p} = \sum_{\dinD} u_{p,d}
\label{B-Eq:f4}
\end{align}
Objective \ref{B-enum:5} is to minimize the number of fractions scheduled on a non-preferred machine stated by the treatment protocol. Therefore, the variable $v_{p,f}$ is introduced and \eqref{B-Eq:IP_machinePref} is used to compute the fractions where a patient is scheduled on a non-preferred machine. The preference violations are summed in \eqref{B-Eq:f5}.
\begin{align}
f_{5,p} = \sum_{\finF_p} v_{p,f}
\label{B-Eq:f5}
\end{align}
There is a cost for switching to a machine that is only a partially beam-matched, which should be minimized according to objective \ref{B-enum:6}. If fraction $f$ is scheduled on a machine in a group of completely beam-matched machines, but $f+1$ is not, then it must be scheduled on a partially matched machine by \eqref{B-Eq:q(p,m,k,f)=q(p,m,k+1,f+1)}. The variable $s_{p,f}$ is one if there is a switch to a partially matched machine, enforced by constraint \eqref{B-Eq:IP_machineswitch}. All machine switches to partially matched machines are summed in \eqref{B-Eq:f6}.
\begin{align}
f_{6,p} = \sum_{f=1}^{F_p-1} s_{p,f}
\label{B-Eq:f6}
\end{align}

\subsection{Column Generation IP Model}
\label{B-Sec:ColGen}
\begin{figure*}[ht]
\small
\begin{alignat}{3}
	 &\text{minimize} \quad        &  & \mathmakebox[0pt][l]{1+\sum_{\pinP} \sum_{i \in \mathcal{K}_p} c_{p,i} a_{p,i}}                                                                                     \label{B-Eq:master obj}\\
	 & \text{subject to} \quad &  &  \sum_{i \in \mathcal{K}_p} a_{p,i} = 1\,, \enskip                                                                                                                                                                      &  & \forall \pinP  \label{B-Eq:master a=1}\\
	 &                   &  & \begin{multlined}\sum_{\pinP} \sum_{i \in \mathcal{K}_p} a_{p,i} \Big((x_{p,m,d,w}^i - t_{p,m,d,w}^i) \mathit{dur}_p  + t_{p,m,d,w}^i \mathit{dur}_{p0}\Big)  \\[-1ex] +S_{m,d,w} \leq L_w,\end{multlined} \enskip                                                                                                                                                            &  &  \begin{multlined} \forall \minM, \dinDw, \winW   \end{multlined}        \label{B-Eq:master sched}      \\
	 &                   &  & \begin{multlined}\sum_{i \in \mathcal{K}_{p}} a_{p,i} \sum_{\dinDw} d \sum_{\minM} \sum_{\winW} t_{p,m,d,w}^i \leq \\ \sum_{i \in \mathcal{K}_{p+1}} a_{p+1,i}\sum_{\dinDw} d \sum_{\minM} \sum_{\winW} t_{p+1,m,d,w}^i  \end{multlined}\enskip & & \begin{multlined}
	 \forall \hinH,\\ p \in \mathcal{P}_h \text{ where } d_{L,p} \leq d_{L,p+1}
	 \end{multlined}  \label{B-Eq:master_dominance}\\	 
	 &                   &  & a_{p,i} \in \{0,1\} \enskip & & \forall \pinP, i \in \mathcal{K}_p                                      \label{B-Eq:master integer}
\end{alignat}
\end{figure*}

The problem is reformulated as a set covering model, where the decision variables represent schedules for each patient. Each patient has an associated index set $\mathcal{K}_p$ of feasible schedules, and the variable $a_{p,i}=1$ if schedule $\iinK$ is allocated to $\pinP$, and 0 otherwise. Since generating all feasible schedules would be too expensive, a column generation approach is presented, which consists of a (restricted) master problem and one subproblem for each patient $\pinP$. The master problem is the schedule selection problem, which is solved to make the overall schedule feasible and optimal. In the subproblems, for each patient a new schedule is generated that fulfills all medical and technical constraints, and the sets of feasible schedules are dynamically updated by the column generation procedure presented in Algorithm~\ref{B-Alg:ColGen}. The algorithm gives a nearly optimal solution, but since the problem is converted from a linear program to an IP in the last step, some schedules not generated by the procedure could potentially improve the integer solutions. The algorithm to generate the initial schedules is presented in Algorithm~\ref{B-Alg:ColGenGreedy}. The number of initial schedules is set to 75, as a larger number does not seem to decrease solution times. 
\begin{algorithm}[ht]
	\caption{Column generation} 
	\begin{algorithmic}[1]
	\State Generate a reduced set of schedules $\mathcal{K}_p'$ for each patient $\pinP$ using Algorithm \ref{B-Alg:ColGenGreedy}
	\While {schedule with negative reduced cost exists}
	    \State Solve linear relaxation of restricted master problem
	    \For {$\pinP$}
	        \State \Longunderstack[l]{Update subproblem objective \eqref{B-Eq:colgen sub objective} with dual\\ variables $\lambda, \gamma$ and $\eta$ from continuous restricted \\master problem solution and solve}         
	        \If {negative reduced cost} 
	            \State \Longunderstack[l]{Add new schedule to $\mathcal{K}_p'$. Coefficients of the\\ new column are given by the optimal\\ solution vector to the subproblem}
	         \EndIf
	    \EndFor
	\EndWhile
	\State Solve restricted master problem using integer values
	\end{algorithmic} 
	\label{B-Alg:ColGen}
\end{algorithm}

\begin{algorithm}
	\caption{Schedule generation procedure} 
	\begin{algorithmic}[1]
	\State Sort $\pinP$ by increasing order of the day limits $d_{L,p}$
	\For {number of initial schedules}
	\For {each protocol $\hinH$ in randomized order}
	    \For {patient $p \in \mathcal{P}_h$ in sorted order}
	        \State Assign $p$ to a random machine $\minM_p$
	        \State \Longunderstack[l]{Schedule $p$ randomly on one of the two first \\available days in a random time window\\ $\winW$ according to the partially occupied\\ input schedule $S$}
	         \State Insert feasible schedule to $\mathcal{K}_p$
	    \EndFor
	\EndFor
	\EndFor
	\end{algorithmic} 
	\label{B-Alg:ColGenGreedy}
\end{algorithm}
\vspace{-10pt}

\mypar{Master problem}
Model \eqref{B-Eq:master obj}-\eqref{B-Eq:master integer} is the master problem: the restricted master
problem is made of a subset $\mathcal{K}_p' \subset \mathcal{K}_p$ of feasible schedules for each $\pinP$. A column in the master problem corresponds to a feasible schedule $\iinK$ for patient $\pinP$. The pure IP variables are now parameters: $x_{p,m,d,w}^i$, $t_{p,m,d,w}^i$, $q_{p,m,d,f}^i$, $y_{p,d,w}^i$, $z_{p,d}^i$, $u_{p,d}^i$, $v_{p,f}^i$, and $s_{p,f}^i$ for $\pinP, \minM, \dinDw, \finF_p, \winW$ have fixed values that satisfy the scheduling constraints presented in Section \ref{B-Sec:IP pure}. Each schedule has an associated (fixed) cost $c_{p,i}$ that is computed using \eqref{B-Eq:f1}-\eqref{B-Eq:f6} with weights $\alpha_1$, $\dots$, $\alpha_6$ according to \eqref{B-Eq:colgen c_pi}:
\begin{alignat}{1}
\begin{multlined}c_{p,i} = \alpha_1 f_{1,p}^i+ \alpha_2 f_{2,p}^i+ \alpha_3 f_{3,p}^i+ \\  \alpha_4f_{4,p}^i + \alpha_5 f_{5,p}^i+ \alpha_6f_{6,p}^i .\end{multlined}
\label{B-Eq:colgen c_pi}
\end{alignat}
The objective function \eqref{B-Eq:master obj} states that the aim is to minimize the total cost of the chosen schedules, plus an offset of 1 to make it equivalent with the other models. Constraint \eqref{B-Eq:master a=1} states that exactly one schedule is chosen for each patient. Constraint \eqref{B-Eq:master sched} ensures that all chosen schedules will fit in the schedule. Constraint \eqref{B-Eq:master_dominance} states that the start day of a patient with shorter day limit should always be before or equal to the start day of a patient with longer day limit if they have the same treatment protocol, by multiplying the master variable with the start day of the corresponding schedule. Note the abuse of notation, where $p+1$ denotes the next entry in $\mathcal{P}_h$.

Relaxing the integer assumption and solving the LP yields the dual variables $\lambda_p$ associated with (\ref{B-Eq:master a=1}), $\gamma_{m,d,w}$ associated with (\ref{B-Eq:master sched}) and $\eta_{h,p}$ associated with \eqref{B-Eq:master_dominance}.

\mypar{Subproblems}
One subproblem is formed for each patient $\pinP$, with the aim to generate a new feasible schedule to add to $\iinK'$, i.e., as a column to the restricted master problem. The constraints are the same as the pure IP formulation \eqref{B-Eq:q(p,m,k,f)=q(p,m,k+1,f+1)}-\eqref{B-Eq:sum(x) = sum(q)}, \eqref{B-Eq:MIP_y_1}-\eqref{B-Eq:MIP_integers}. The schedule availability constraint \eqref{B-Eq:IP x*dur+S<=L_w} is replaced by \eqref{B-Eq:colgen sub sched}, since the subproblems only deal with one patient at a time. 
\begin{alignat}{2}
\begin{multlined}(x_{p,m,d,w}^i - t_{p,m,d,w}^i) \mathit{dur}_p +\\ t_{p,m,d,w}^i \mathit{dur}_{p0}+  S_{m,d,w} \leq L_w, \end{multlined} \enskip                                                                                                              &  \begin{aligned}\forall &\minM, \\ &\dinDw,\\&\winW \end{aligned} \label{B-Eq:colgen sub sched}
\end{alignat}

The subproblem objective function \eqref{B-Eq:colgen sub objective} is the cost of the schedule defined by \eqref{B-Eq:colgen c_pi}, minus the master dual variables $\lambda_p$, $\gamma_{mdw}$ and $\eta_{h,p}$ multiplied by the coefficients given from their respective constraints in the master problem. The dual variable $\eta_{h,p}$ associated with \eqref{B-Eq:master_dominance}, with $h$ being the protocol of patient $p$, is partly shared between patient $p$ and $p+1$ for $p \in \mathcal{P}_h$, because of the formulation of \eqref{B-Eq:master_dominance} (using the same abuse of notation).
\begin{small}
\begin{align}
\begin{multlined}\text{minimize} \quad c_{p,i} -  \lambda_p - \sum_{\minM} \sum_{\dinDw} \sum_{\winW} \gamma_{m,d,w} \Big(\\(x_{p,m,d,w}^i - t_{p,m,d,w}^i) \mathit{dur}_p 
+  t_{p,m,d,w}^i \mathit{dur}_{p0}\Big)  - \\ (\eta_{h,p-1}-\eta_{h,p})\sum_{\minM} \sum_{\dinDw} \sum_{\winW} t_{p,m,d,w}^i d \end{multlined} \label{B-Eq:colgen sub objective}
\end{align}
\end{small}

Since the subproblems are isolated from each other, constraint \eqref{B-Eq:colgen sub sched} can be satisfied as long as the input schedule $S_{m,d,w}$ together with the current patient's duration do not require more than the entire window capacity $L_w$. From this point of view, the subproblems are very easy to solve. On the other hand, the optimization is exclusively guided by the values of the dual variables which might lead to a larger number of iterations.

\subsection{Constraint Programming Model}
\label{B-Sec:CPmodel}
\begin{figure*}[b]
\small
\begin{alignat}{3}
& \text{minimize} \quad   && \mathmakebox[0pt][l]{
1+ \sum_{\pinP}( \alpha_1 k_{1,p} + \alpha_2 k_{2,p} + \alpha_3 k_{3,p} + \alpha_4 k_{4,p} + \alpha_5 k_{5,p} + \alpha_6 k_{6,p} )} \label{B-Eq: CP Objective}\\
& \text{subject to} \quad&& \sd{p} \in \mathcal{A}_p \enskip && \forall \pinP
	\label{B-Eq: CP allowed start day}\\
&                         && \machine{p}{d} \in \mathcal{M}_p \enskip && \forall \pinP, \dinD_w 
	\label{B-Eq: CP machine} \\
&                         && \mgP = b_\mathcal{M} \text{ if } \machine{p}{\sd{p}} \in b_{\mathcal{M}} \enskip && \forall \pinP, \binBm
	\label{B-Eq: CP machine group} \\
&                         && \machine{p}{d} \in \mgP \enskip && \forall \pinP, \dinDw
	\label{B-Eq: CP machine group2} \\	
&                         && d_{\min,p} \leq \sd{p} \leq D_w - F_p + 1 \enskip && \forall \pinP
	\label{B-Eq: CP_daymin} \\
&                         && \begin{aligned} & \window{p}{d} > 0 \iff 
                            d \geq \sd{p}  \land
                            d < \sd{p} + F_p \end{aligned} \enskip &&  \begin{aligned} \forall \pinP, \dinDw \end{aligned}
	\label{B-Eq: CP window zero} \\
&                         && \machine{p}{d} > 0 \iff \window{p}{d} > 0 \enskip && \forall \pinP, \dinDw
	\label{B-Eq: CP machine_window} \\	
&                         && \sum_{\dinDw} \machine{p}{d} = F_p \enskip && \forall \pinP
	\label{B-Eq: CP number_fractions} \\	
&                         &&\begin{aligned}
\texttt{bin}\texttt{\_packing}(&[\infty, L_1, \dots, L_W], \\
                      &[\text{for each } \pinP: \\
                      &\text{if } m=\machine{p}{d} \land d=\sd{p} \text{ then } \\ 
                      & \quad (\window{p}{d}, \textit{dur}_{p0} ) \\
                      &\text{else if } m=\machine{p}{d} \text{ then } \\
                      & \quad (\window{p}{d}, \textit{dur}_{p} ) \\
                      &\text{else } (0,0)] ++ \\
                      &[\text{for each } \winW: (w,S_{m,d,w})]) \end{aligned} \enskip && \begin{aligned} &\forall \minM, \dinDw \\ & \\ & \\& \\&\end{aligned} \label{B-Eq: Cp bin packing}\\
&                     && \begin{multlined}\sd{p} \leq  \sd{p+1} \\ \end{multlined}
	\enskip && \begin{multlined}\forall	 h \in \mathcal{H}, \pinP_h \text{ where } d_{L,p} \leq d_{L,p+1}\end{multlined}
	\label{B-Eq: CP dominance}
\end{alignat}
\end{figure*}
\normalsize

In \cite{Frimodig2019}, the authors found that the CP model that used bin packing constraints was more efficient than the CP model that used scheduling constraints for the RT patient scheduling problem. Therefore, we present a CP model that uses a global bin packing constraint. The variables in the CP model are presented in Table \ref{B-tab:variables_CP}. The aim is to assign a $\sd{p}$, a $\machine{p}{d}$, and a time $\window{p}{d}$ to each patient $\pinP$ for each day $\dinDw$ using the formulation \eqref{B-Eq: CP Objective}-\eqref{B-Eq: CP dominance}. The variable $\mgP$ is a function of $\machine{p}{\sd{p}}$.

\begin{table}
\caption{Variables in the CP model}
\label{B-tab:variables_CP}
\begin{tabular}{P{0.42\linewidth} P{0.48\linewidth}}
\hline\noalign{\smallskip}
    	$\window{p}{d} \in \{0,\dots, W\}$ & Window patient $\pinP$ is scheduled in on day $\dinDw$, where 0 represents patient $p$ not being scheduled day $d$ \\
	$\machine{p}{d} \in \{0,\dots, M\}$    & Machine patient $\pinP$ is scheduled on day $\dinDw$, where 0 represents patient $p$ not being scheduled day $d$  \\
	$\sd{p} \in \Dw$                         & Start day for patient $\pinP$ \\
	$\mgP \in \mathcal{B_M}$        & Group of beam-matched machines patient $\pinP$ is scheduled on\\
\noalign{\smallskip}\hline
\end{tabular}
\end{table}

\mypar{Constraints}
Constraint \eqref{B-Eq: CP allowed start day} states that each patient must start on an allowed start day. Constraint \eqref{B-Eq: CP machine} states that the treatment must be scheduled on a machine allowed for that patient. Constraint \eqref{B-Eq: CP machine group} and \eqref{B-Eq: CP machine group2} make all fractions scheduled on machines from a group of beam-matched machines. Constraint \eqref{B-Eq: CP_daymin} limits the earliest start day and that the start day for each patient should be at least $F_p$ days from the end of the planning horizon. Equation \eqref{B-Eq: CP window zero} states that the variable \window{p}{d} should be nonzero if and only if treatment has started but not ended. The active machine days should be the same as the active window days, which is stated in \eqref{B-Eq: CP machine_window}. The number of active machine days should be the same as the number of fractions and is stated in \eqref{B-Eq: CP number_fractions}. This is a redundant constraint, as it is already enforced by \eqref{B-Eq: CP window zero} and \eqref{B-Eq: CP machine_window}, but added as it helps performance during search. 

To ensure that the patients fit in each window, a global bin packing constraint \cite{Shaw2004} is used. In \eqref{B-Eq: Cp bin packing}, the first line states that the capacity of window~0 is infinite (corresponding to not being scheduled). Window $1,\dots,W$ have capacity $L_1,\dots,L_W$. In \eqref{B-Eq: Cp bin packing}, the bin choice and required allocation are created together as a list of pairs $(bin,size)$ for each patient. The first value in the pair corresponds to the bin, here the $\window{p}{d}$, and the second is the size of the item, which corresponds to the duration of the treatment. If the patient is not scheduled on the particular machine, window~0 is chosen with item size~0. This list is concatenated (using the $++$ operator) with a list of pairs used to include already occupied timeslots $S_{m,d,w}$ in window $w$. 

The same dominance breaking as in the IP model is included: if two patients have the same treatment protocol, the one with shorter treatment target should start treatment first, which is stated in Constraint \eqref{B-Eq: CP dominance}, using the same abuse of notation as before, i.e., $p+1$ denotes the next entry in $\mathcal{P}_h$.

\mypar{Objective function}
Equation (\ref{B-Eq: CP Objective}) shows the generalized objective function, which is divided into six parts $k_1,\dots,k_6$ according to Section \ref{B-Sec: Problem Formulation}, and combined with weights $\alpha_1,\dots,\alpha_6$. An offset of 1 is included to make it equivalent to the IP objective function \eqref{B-Eq:IP objective}.

Objective \ref{B-enum:i} is to minimize a weighted sum of the waiting times, which is done in \eqref{B-Eq:f1_CP} by penalizing the number of days between the first allowed start day $d_{\min,p}$ and the start day, multiplied with weight $c_p$ corresponding to the priority group of patient $p$.  
\begin{align}
 k_{1,p} = c_p ( \sd{p} - d_{\min,p})
 \label{B-Eq:f1_CP}
\end{align}
Objective \ref{B-enum:ii}, to minimize a weighted sum of the violations of the waiting time targets, is formulated in \eqref{B-Eq:f2_CP}. The target violation is zero if the start day is before the waiting time target, and otherwise penalized linearly. 
\begin{align}
 k_{2,p} = c_p \max(0, \sd{p} - d_{L,p})
 \label{B-Eq:f2_CP}
\end{align}
Objective \ref{B-enum:iii} is to minimize the number of window switches. Therefore, a penalty of value one is added each time the window is switched. Since only the days when treatment has started but not finished are relevant, i.e., when $\window{p}{d} \neq 0$, we let the active treatment days form a set $\mathcal{D}_a \subset \Dw$ and compute the number of window switches on that set in \eqref{B-Eq:f3_CP}. Note the abuse of notation, where $d+1$ denotes the next entry in $\mathcal{D}_a$.
\begin{align}
 k_{3,p} = \sum_{d \in \mathcal{D}_a} (\window{p}{d} \neq \window{p}{d+1})
 \label{B-Eq:f3_CP}
\end{align}
Objective \ref{B-enum:iv} is to minimize the violations of the window preferences, which is formulated in \eqref{B-Eq:f4_CP} for $\pinPpref$ (otherwise $k_{4,p}=0$).
\begin{align}
 k_{4,p} = \sum_{d \in \mathcal{D}_a} |\window{p}{d}-w^{\mathit{pref}}_{p}|
 \label{B-Eq:f4_CP}
\end{align}
Objective \ref{B-enum:5}, to minimize the fractions scheduled on a non-preferred machine, is formulated in \eqref{B-Eq:f5_CP}.
\begin{align}
 k_{5,p} =  \sum_{d \in \mathcal{D}_a} (\machine{p}{d} \notin \MprefP)
 \label{B-Eq:f5_CP}
\end{align}
Objective \ref{B-enum:6}, to minimize the number of switches to a machine that is only partially beam-matched, is formulated in \eqref{B-Eq:f6_CP}. It states that if $\machine{p}{d+1}$ is in the set of partially beam-matched machines $p_\mathcal{M}$ for the next day machine $\machine{p}{d}$, then there has been a switch between day $d$ and $d+1$.
\begin{align}
 k_{6,p} = \sum_{d \in \mathcal{D}_a} (\machine{p}{d+1} \in p_\mathcal{M}^{\machine{p}{d}})
 \label{B-Eq:f6_CP}
\end{align}

\subsubsection{Search Heuristic}
In CP, the solvers rely on backtracking algorithms that are used in the tree search-based heuristics. When using backtracking search, a sequence of decisions are made regarding what \textit{variable} to branch on next, and which \textit{value} to assign to the variable. It is well known that the choice of variable and value ordering, also called search heuristic, can be crucial to solving a problem efficiently, see e.g.~\cite{Haralick1980} and \cite{Ginsberg1990}. 

For the CP model in this paper, many different choices of variable and value orderings were investigated. Our initial experiments showed that randomization and restarts are necessary to obtain good results: the restart search helps avoid getting stuck in a non-productive area of the search tree. Several restart strategies were evaluated, and the strategy with best performance was the Luby restart strategy~\cite{Luby}, which gives a specific scheme for when search is restarted.

When deciding what variables to assign random values in the search heuristic, the best performing search heuristic was shown to be tightly related to the objective function \eqref{B-Eq: CP Objective}. It is obvious that assigning the $\sd{p}$ variable a random value will not result in good quality solution, whereas assigning $\machine{p}{d}$ a random value will give the benefit of a wider search tree. 
The search heuristic is described in Algorithm \ref{B-Alg:CP_Search}. Large Neighborhood Search \cite{Shaw1998} was also tested but did not improve overall performance.
\vspace{-2pt}
\begin{algorithm}
	\caption{CP Search Heuristic} 
	\begin{algorithmic}[1]
	\State The Luby restart strategy \cite{Luby} with parameter 75 is used
	\State Assign the sum $\sum_{\pinP}$ of all window switches $k_{3,p}$  \eqref{B-Eq:f3_CP} to the minimum value
	\State Assign the sum $\sum_{\pinP}$ of fractions on non-preferred machine $k_{5,p}$ \eqref{B-Eq:f5_CP} to the minimum value
	\State Assign the sum $\sum_{\pinP}$ of machine switches to partially beam-matched machines $k_{6,p}$  \eqref{B-Eq:f6_CP} to the minimum value
	\State Create patient list sorted by priority group, with all priority A patients first, followed by priority B and C
	\For {$p$ in sorted list of patients}
        \State Assign $\sd{p}$ to the minimum value
        \State Assign $\machine{p}{d}$ a random value for each $\dinDw$
        \State \Longunderstack[l]{Assign the window preference violation $k_{4,p}$ in \eqref{B-Eq:f4_CP} \\ to the minimum value}
        \vspace{2pt}
        \State Assign $\window{p}{d}$ a random value for each $\dinD$
	\EndFor
	\end{algorithmic} 
	\label{B-Alg:CP_Search}
\end{algorithm}

\vspace{-15pt}
\subsection{Combined CP/IP method}
\label{B-Sec:Combined}

For difficult MIP problems, providing the solver a good quality input solution ("warm start") can improve performance. The solver processes the input solution before starting branch-and-cut to get a lower/upper bound to use during the optimization, which allows it to eliminate parts of the search space. For this reason, the CP model is used to find a feasible solution to use as a warm start in the IP model. The first feasible solution found by the CP model is generally of good quality because of the search heuristic, which gives the MIP solver a useful upper bound during branch-and-cut. The CP solution is transformed to the format of the IP $x$-variables by letting $x_{p,m,d,w}=1$ if and only if $w=\window{p}{d}$ and $m=\machine{p}{d}$ for $\pinP, \minM, \dinDw, \winW$. This is provided to the MIP-solver using the built-in functionality for advance starting.

\section{Experimental setup}
\label{B-Sec: Experiments}
This section presents the setup for the experiments. Section \ref{B-sec:time horizon} presents how the time horizon is computed. Section \ref{B-sec:data_benchmarks} describes the historical patient data from Iridium Netwerk and how the problem instances are generated. Section \ref{B-sec:objective_functions} presents the objective functions.

The experiments are run on a Windows 10 machine with an Intel\textsuperscript{\textregistered} Core{\texttrademark} i9-7940X X-series processor and 64~GB of RAM. The patient arrival model used when generating benchmarks is built with Python~3.8. The IP models are solved using the MIP solver of CPLEX~12.10 in the Python~API with the default parameters.
The CP model is written in MiniZinc~2.5.5 \cite{MiniZinc}, and uses the Gecode~6.3.0 solver \cite{Gecode}. Other CP solvers were tested, such as the lazy clause generation solver Chuffed \cite{Chuffed}, but Gecode gave the best overall results on the tested problem instances. 
The maximum allowed CPU time was set to 1 hour per run.

\subsection{Computing the time horizon $D_w$}
\label{B-sec:time horizon}
The models all depend on the number of days in the time horizon. $D_w$ should be large enough to schedule all treatments, but a larger $D_w$ may weaken performance due to larger problem dimensions. A heuristic to compute $D_w$ is presented in Algorithm~\ref{B-Alg:Dw}. A random schedule is computed, and the value of $D_w$ is set to the last utilized day out of all patients, plus 30 days that are added to augment the search space.

\begin{algorithm}
	\caption{Time Horizon Heuristic} 
	\begin{algorithmic}[1]
	\State Sort patients by their priority
	\State Sort patients of same priority by the day of arrival
		\For {$p$ in sorted list of patients}
		    \State Randomly assign $p$ to a machine in $\mathcal{M}_p$
		    \State \Longunderstack[l]{Schedule $p$ randomly on one of the two first available \\days in the first available time window $\winW$ \\according to the partially occupied input schedule $S$}
		    \State $D_p=$ the last day in the schedule 
		\EndFor
	\State $D_w = \max_p(D_p) + 30$    \Comment{Augment search space by 30}
	\end{algorithmic} 
	\label{B-Alg:Dw}
\end{algorithm}

\subsection{Generating problem instances using historical clinic data}
\label{B-sec:data_benchmarks}

In 2020, 4070 patients received 5500 treatments at Iridium Netwerk. From the historical data, the empirical distribution of the 72 different treatment protocols can be computed. Each protocol states the machines that are equipped for treating the particular tumor, what machines that are preferred for treating the target, and the duration for the first and subsequent fractions. The average number of fractions for each treatment protocol is computed from historical data. Furthermore, each treatment protocol is given a priority (A, B or C) by a radiation oncologist (MD) at Iridium, which will give an equivalent patient priority. In 2020, Iridium Netwerk operated 10 linacs and 255 days were used to treat patients, resulting in an average arrival rate of 16 patients per working day. No records were kept over the patients' time window preferences, but the booking administrators estimate that $80\%$ of the patients have a preference, of which $65\%$ prefer a treatment before noon and $35\%$ prefer the afternoon, and that $20\%$ of the patients have no preference.

Literature shows the majority of patients find it reasonable to receive a notification of the treatment three days in advance \cite{Olivotto2015}. Therefore, in this paper the duration of notice is three days for priority~B and C patients, while priority~A patients are notified immediately. All fractions are communicated and cannot be re-planned, as this is the current practice at Iridium Netwerk. The schedule can change until being communicated; booking decisions are postponed to the next day for patients scheduled after the notification period. The notification period length is straightforward to change.

\begin{figure*}
    \centering
    \includegraphics[width=17cm]{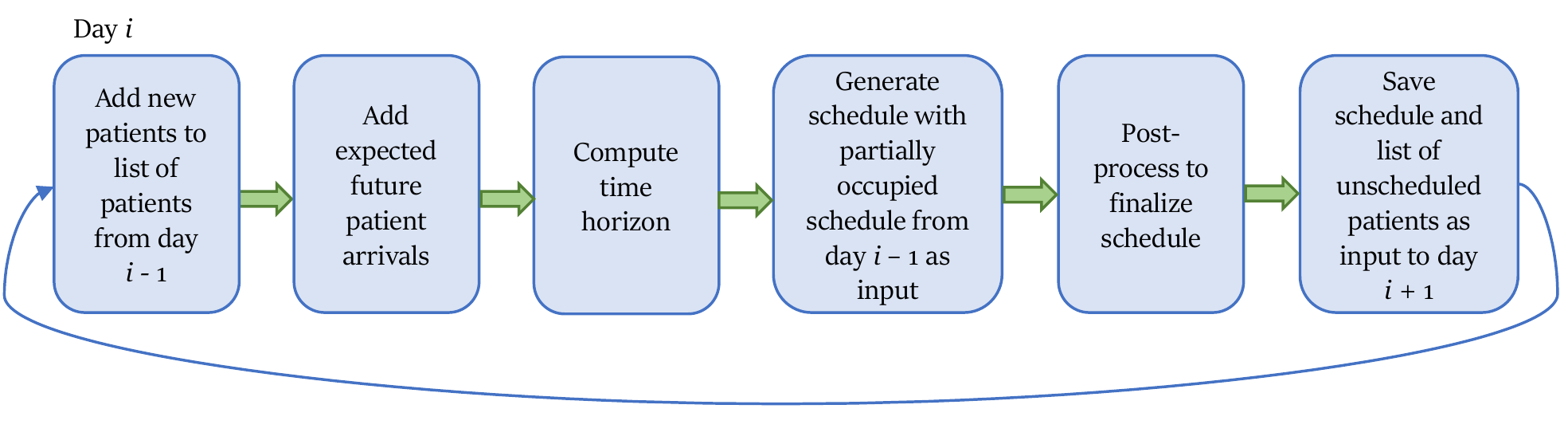}
    \caption{Problem generation algorithm. Each day $i$, the patient arrivals are simulated and a schedule is computed. The resulting schedule and list of unscheduled patients are saved as input to day $i+1$}
    \label{B-fig:PAM}
\end{figure*}

\mypar{Problem generation} A model for patient arrivals is developed. The goal is to mimic scheduling behavior to generate realistic problem instances based on the historical data from Iridium Netwerk. Each problem instance should represent different scenarios, altering the number of patients to be scheduled and the partially occupied input schedule. The problem generation algorithm simulates each day at the clinic; the patients arriving, and the resulting schedules. An overview of the steps during each simulated day can be seen in Figure \ref{B-fig:PAM}. 

In the first step, new patients are assumed to arrive according to a Poisson process based on historical arrival rates. Each patient is randomly assigned a treatment protocol from the empirical distribution of protocols. Secondly, priority~A and B patients that are expected to arrive in the coming four weeks are added to the problem as placeholder (\textit{dummy}) patients. Their treatment target dates and earliest start day are set from when they are expected to arrive. Thirdly, Algorithm \ref{B-Alg:Dw} is run to determine the time horizon. Next, the IP model is run to generate a schedule (but any model could be used in the problem generation phase). The schedule assigns each patient a machine, treatment days, and time windows. Next, the results are post-processed. 
Patients that start treatment within the duration of notice are fixed to the schedule, whereas the booking decisions are postponed to the next day for patients that are scheduled more than three days away. Finally, the schedule is saved and is given together with the list of unscheduled patients as input to the next day.

\mypar{Problem benchmarks} 
In order to evaluate how well the models scale to different problem sizes, the problem instance generator is used for four setups: two average arrival rates, $\lambda = \{16,18\}$, and two different number of time windows, $W=\{2,4\}$. When $W=2$, the time window preferences from Iridium Netwerk are used. For $W=4$, an estimation is that $25\%$ prefer the first window, $25\%$ prefer the last window, and $50\%$ have no preference. For each setup, 20 days are randomly chosen between day 50 and 300 in the simulation to form the problem benchmarks. These instances represent different scenarios, altering the patient flow and the partially occupied input schedule. 

In the generated problem benchmarks, the number of patients to schedule, including expected future arrivals as placeholder  patients, varies. When $\lambda=16$,  $P \in [225,239]$ with an average of $230.2$, and when $\lambda=18$, $P \in [250,268]$ with an average of $256.5$. The time horizon $D_w$ also varies; $D_w \in [79,89]$ when $\lambda=16$, and $D_w \in [78,94]$ when $\lambda=18$. The average occupancy on all machines except M10 (which is specialized and always has a lower occupancy) of the first day is $65.7\%$ for $\lambda=16$ and $73.5\%$ when $\lambda=18$. 
All problem instances used in this paper are publicly available\footnote{Access through this link:  \url{https://osf.io/45qw2/?view_only=4a0a67e21cb542df8f9a0f74241de825}}.

\subsection{Objective functions}
\label{B-sec:objective_functions}
 As presented in Section \ref{B-Sec: Problem Formulation}, there are several objectives considered: \ref{B-enum:i} is to minimize a weighted sum of the waiting times, \ref{B-enum:ii} is to minimize a weighted sum of the violations of the target dates, \ref{B-enum:iii} is to minimize the number of time window switches, \ref{B-enum:iv} is to minimize violations of time window preferences, \ref{B-enum:5} is to minimize the number of fractions scheduled on non-preferred machines, and \ref{B-enum:6} is to minimize the number of times a patient switches between machines that are only partially beam-matched. These objectives are combined into different objective functions using weights $\alpha_1, \dots,\alpha_6$ as presented in \eqref{B-Eq:IP objective}, \eqref{B-Eq:colgen c_pi} and \eqref{B-Eq: CP Objective}. The different combinations that form the objective functions used in the experiments are presented in Table \ref{B-tab:obj_funs}. 
 
 \begin{table}
\caption{The different objective function combinations}
\label{B-tab:obj_funs}
\begin{tabular}{ P{0.152\linewidth} P{0.247\linewidth} P{0.452\linewidth} }
\hline\noalign{\smallskip}
 Objective function number & Combination & Weights \\
\noalign{\smallskip}\hline\noalign{\smallskip}
\#1 & \ref{B-enum:i}+\ref{B-enum:ii}+\ref{B-enum:iii}+ \ref{B-enum:5}+\ref{B-enum:6}  & $\alpha_1=50$, $\alpha_2=100$, $\alpha_3=1$, $\alpha_4=0$, $\alpha_5=10$, $\alpha_6=10$  \\
\#2&  \ref{B-enum:i}+\ref{B-enum:ii}+\ref{B-enum:iii}+ \ref{B-enum:iv} &$\alpha_1=50$, $\alpha_2=100$, $\alpha_3=1$, $\alpha_4=1$, $\alpha_5=0$, $\alpha_6=0$\\
\#3&  \ref{B-enum:i}+\ref{B-enum:iii}+\ref{B-enum:5} &$\alpha_1=100$, $\alpha_2=0$, $\alpha_3=1$, $\alpha_4=0$, $\alpha_5=10$, $\alpha_6=0$\\
\#4&  \ref{B-enum:i}+\ref{B-enum:iii}+\ref{B-enum:iv}+ \ref{B-enum:5}+\ref{B-enum:6} &$\alpha_1=100$, $\alpha_2=0$, $\alpha_3=1$, $\alpha_4=5$, $\alpha_5=10$, $\alpha_6=10$\\
\noalign{\smallskip}\hline
\end{tabular}
\end{table}
The objective functions are designed to mimic the scheduling policies at different clinics. For example, in some countries there are no official treatment target dates, and therefore objective \ref{B-enum:ii} is not active in objective functions \#3 and \#4. Some clinics do not consider the patients' preferences of treatment time of the day, which is why objective \ref{B-enum:iv} is not included in \#1 and \#3. For clinics that do not have multiple hospitals, it is unlikely that the problem with machine switches to partially beam-matched machines exists, hence objective \ref{B-enum:6} is irrelevant and therefore not included in \#2 and \#3. Some clinics may not state preferred machines, thus \ref{B-enum:5} is not included in objective function \#2.

The waiting time has a large negative effect on the patient outcome, especially for acute patients, see e.g. \cite{CHEN2008}. Therefore, both objective \ref{B-enum:i} and \ref{B-enum:ii} also have the weight $c_p$ for each patient (see \eqref{B-Eq:f1} and \eqref{B-Eq:f2}), which reflects the severeness of delaying treatment start for the different priority groups.  In objective functions \#1 and \#2, the weights of $\alpha_1$ and $\alpha_2$ show that if the patients are of the same priority group, it is never desirable to minimize a patient's waiting time at cost of another patient missing their treatment target date. However, if one patient is priority A and one is priority C, the latter is allowed to violate the treatment target date if it means the priority A patient gets a shorter waiting time. Furthermore, the weights of $\alpha_5$ and $\alpha_6$ compared to $\alpha_2$ in objective function \#1 indicate that it is preferred for a patient to start their treatment earlier at the cost of either switching linacs or scheduling on non-preferred linacs. Finally, the weight of $\alpha_3$ is lower than the rest; if possible, all fractions should be scheduled in the same time window, but never at the cost of any of the other objectives.

Objective function number \#4 is most similar to what is used at Iridium Netwerk today. There are no official treatment target dates in Belgium, therefore objective \ref{B-enum:ii} is not active. Minimizing waiting times is by far the most important objective, thus the weight of this objective, $c_p\alpha_1$, is the largest. It is more important to fulfill the patient preferences regarding time windows than to schedule them in the same time window every day, thus $\alpha_3 < \alpha_4$. At Iridium, the current practice is to never schedule patients with switches between partially beam-matched machines. However, the staff at Iridium agrees that this should be allowed if it will lead to a minimized waiting time. Therefore, the penalty for scheduling patients on a non-preferred machine is the same as for scheduling patients with machine switches between partially beam-matched machines. If a treatment starts on a non-preferred machine, the aim is to switch to a preferred machine as soon as possible. This is true also in \#4: the cost of switching to a partially beam-matched preferred machine will be lower when there is more than one fraction left to schedule, since the switch is a one-time cost.

\section{Results}
\label{B-Sec: Results}
The four different setups ($\lambda = \{16,18\}$, $W=\{2,4\}$) are run with the four objective functions presented in Table \ref{B-tab:obj_funs}, giving a total of 16 different combinations that have 20 problem instances each.

\begin{table*}
\caption{Computational time results. For each combination of arrival rate, number of time windows, objective function and model, column \textbf{A} shows the average CPU time, column \textbf{B} presents the median CPU time, and column \textbf{C} shows the cumulative CPU time for the 20 instances. The timeout is set to 1 hour of CPU time} 
\label{B-tab:results-time}
\setlength\extrarowheight{2pt}
\begin{tabular}{P{0.055\linewidth}|P{0.068\linewidth}| P{0.077\linewidth} | P{0.036\linewidth} P{0.036\linewidth} P{0.036\linewidth} | P{0.036\linewidth} P{0.036\linewidth} P{0.036\linewidth} | P{0.036\linewidth} P{0.036\linewidth} P{0.036\linewidth} |P{0.036\linewidth} P{0.036\linewidth} P{0.036\linewidth}}
\hline
Arrival rate &Number of time windows&  Objective function number &\multicolumn{12}{P{0.58\linewidth}}{\textbf{A}: Average CPU time, \textbf{B}: Median CPU time, \textbf{C}: Cumulative CPU time (all in minutes)}  \\
&&& \multicolumn{3}{P{0.1\linewidth}|}{IP} & \multicolumn{3}{P{0.1\linewidth}|}{CP} & \multicolumn{3}{P{0.1\linewidth}|}{CG-IP} & \multicolumn{3}{P{0.15\linewidth}}{Combined CP/IP}  \\
&&& \textbf{A} & \textbf{B} & \textbf{C} & \textbf{A} & \textbf{B} & \textbf{C}  & \textbf{A} & \textbf{B} & \textbf{C} & \textbf{A} & \textbf{B} & \textbf{C}  \\
\hline
\multirow{8}{*}{$\lambda = 16$} &\multirow{4}{*}{$W=2$}
&   \#1 & 7 & 7 & 155 &  21 & 11 & 435 &  2 & 2 & 51 &    12 & 7 & 247 \\
&& \#2 & 2 & 2 & 57 &  51 & 60 & 1039 &  2 & 2 & 53 &    12 & 5 & 249 \\
&& \#3 & 2 & 3 & 56 &  30 & 16 & 600 &  2 & 2 & 47 &    14 & 4 & 285 \\
&& \#4 & 5 & 5 & 113 &  52 & 60 & 1046 &  3 & 3 & 63 &   14 & 12 & 289 \\
\cline{2-15}
&\multirow{4}{*}{$W=4$}
&   \#1 & 42 & 40 & 856 &  24 & 15 & 482 &  3 & 3 & 70 &   33 & 38 & 672 \\
&& \#2 & 6 & 6 & 136 &  60 & 60 & 1200 &  8 & 6 & 163 &   28 & 36 & 560 \\
&& \#3 & 5 & 5 & 111 &  42 & 40 & 843 &  3 & 3 & 78 &     8 & 7 & 173 \\
&& \#4 & 45 & 53 & 900 &  60 & 60 & 1200 &  9 & 8 & 182 &   47 & 50 & 955 \\
\hline
\multirow{8}{*}{$\lambda = 18$} &\multirow{4}{*}{$W=2$}
&   \#1 & 17 & 13 & 359 &  36 & 34 & 724 &  4 & 3 & 81 &   33 & 41 & 661 \\
&& \#2 & 6 & 3 & 121 &  59 & 60 & 1192 &  3 & 2 & 78 &    18 & 6 & 375 \\
&& \#3 & 7 & 3 & 143 &  45 & 46 & 910 &  6 & 3 & 128 &   20 & 17 & 405 \\
&& \#4 & 14 & 8 & 287 &  57 & 60 & 1149 &  4 & 3 & 91 &   27 & 35 & 547 \\
\cline{2-15}
&\multirow{4}{*}{$W=4$}
&   \#1 & 60 & 60 & 1200 &  40 & 41 & 816 &  8 & 7 & 165 &   49 & 60 & 985 \\
&& \#2 &28 & 25 & 570 &  60 & 60 & 1200 &  48 & 60 & 963 &   40 & 43 & 818 \\
&& \#3 & 9 & 6 & 192 &  54 & 60 & 1093 &  10 & 6 & 204 &    10 & 8 & 207 \\
&& \#4 & 60 & 60 & 1200 &  60 & 60 & 1200 &  55 & 60 & 1114 &  60 & 60 & 1200 \\
\hline
\end{tabular}
\end{table*}

\begin{table*}
\caption{Solution quality results. For each combination of arrival rate, number of time windows, objective function and model, column \textbf{A} shows the mean relative optimality gap, column \textbf{B} presents the median relative optimality gap, and column \textbf{C} shows the proportion of instances that did not have a feasible incumbent solution at timeout for the 20 instances. The timeout is set to 1 hour of CPU time} 
\label{B-tab:results-quality}
\setlength\extrarowheight{2pt}
\begin{tabular}{P{0.055\linewidth}|P{0.068\linewidth}| P{0.077\linewidth} | P{0.036\linewidth} P{0.036\linewidth} P{0.036\linewidth} | P{0.036\linewidth} P{0.036\linewidth} P{0.036\linewidth} | P{0.036\linewidth} P{0.036\linewidth} P{0.036\linewidth} |P{0.036\linewidth} P{0.036\linewidth} P{0.036\linewidth}}
\hline
Arrival rate &Number of time windows&  Objective function number &\multicolumn{12}{P{0.65\linewidth}}{\textbf{A}: Mean relative optimality gap at timeout (\%), \textbf{B}: Median relative optimality gap at timeout (\%), \textbf{C}: Proportion with no feasible solution at timeout (\%)}  \\
&&& \multicolumn{3}{P{0.1\linewidth}|}{IP} & \multicolumn{3}{P{0.1\linewidth}|}{CP} & \multicolumn{3}{P{0.1\linewidth}|}{CG-IP} & \multicolumn{3}{P{0.15\linewidth}}{Combined CP/IP}  \\
&&& \textbf{A} & \textbf{B} & \textbf{C} & \textbf{A} & \textbf{B} & \textbf{C}  & \textbf{A} & \textbf{B} & \textbf{C} & \textbf{A} & \textbf{B} & \textbf{C}  \\
\hline
\multirow{8}{*}{$\lambda = 16$} &\multirow{4}{*}{$W=2$}
&   \#1 & 0.0 & 0.0 & 0 &  0.1 & 0.0 & 0 &  0.0 & 0.0 & 0 &   0.0 & 0.0 & 0 \\
&& \#2 & 0.0 & 0.0 & 0 &  0.9 & 0.0 & 0 &  0.0 & 0.0 & 0 &   0.0 & 0.0 & 0 \\
&& \#3 & 0.0 & 0.0 & 0 &  0.1 & 0.0 & 0 &  0.0 & 0.0 & 0 &   0.0 & 0.0 & 0 \\
&& \#4 & 0.0 & 0.0 & 0 &  1.2 & 0.2 & 0 &  0.0 & 0.0 & 0 &   0.0 & 0.0 & 0 \\
\cline{2-15}
&\multirow{4}{*}{$W=4$}
&   \#1 & 0.2 & 0.0 & 0 &  0.7 & 0.0 & 0 &  0.0 & 0.0 & 0 &   0.0 & 0.0 & 0 \\
&& \#2 &  0.0 & 0.0 & 0 &  8.0 & 7.5 & 5 &  0.0 & 0.0 & 0 &   0.0 & 0.0 & 0 \\
&& \#3 & 0.0 & 0.0 & 0 &  0.7 & 0.0 & 0 &  0.0 & 0.0 & 0 &   0.0 & 0.0 & 0 \\ 
&& \#4 & 0.3 & 0.0 & 30 &  16.5 & 14.6 & 0 &  0.0 & 0.0 & 0 &   0.3 & 0.0 & 5 \\
\hline
\multirow{8}{*}{$\lambda = 18$} &\multirow{4}{*}{$W=2$}
&   \#1 & 0.0 & 0.0 & 10 &  0.5 & 0.0 & 0 &  0.0 & 0.0 & 0 &   0.5 & 0.0 & 5 \\
&& \#2 & 0.0 & 0.0 & 0 &  2.0 & 0.9 & 30 &  0.0 & 0.0 & 0 &   0.0 & 0.0 & 0 \\
&& \#3 & 0.0 & 0.0 & 0 &  0.8 & 0.0 & 20 &  0.0 & 0.0 & 0 &   0.1 & 0.0 & 0 \\
&& \#4 & 0.0 & 0.0 & 10 &  6.1 & 6.4 & 0 &  0.0 & 0.0 & 0 &   0.0 & 0.0 & 5 \\
\cline{2-15}
&\multirow{4}{*}{$W=4$}
&   \#1 & 46.4 & 46.4 & 95 &  4.1 & 0.0 & 0 &  0.0 & 0.0 & 0 &   6.8 & 3.9 & 0 \\
&& \#2 & 0.0 & 0.0 & 0 &  17.4 & 15.0 & 25 &  0.4 & 0.2 & 0 &   0.0 & 0.0 & 0 \\
&& \#3 & 0.0 & 0.0 & 0 &  5.5 & 2.8 & 0 &  0.0 & 0.0 & 0 &   0.0 & 0.0 & 0 \\
&& \#4 & - & - & 100 &  40.0 & 44.3 & 30 &  0.2 & 0.0 & 0 &  44.2 & 46.7 & 40 \\
\hline
\end{tabular}
\end{table*}

\subsection{Computational efficiency}
The models in Section \ref{B-Sec: Models} are run for the 20 problem instances for each of the 16 combinations. Table~\ref{B-tab:results-time} presents the mean, median and cumulative CPU times. Table~\ref{B-tab:results-quality} presents the quality of the solutions; it shows the mean and median of the relative optimality gap, i.e., the mean or median of $(x-y)/y$, where $x$ is the current best objective value and $y$ is the proven optimal value. The table also presents the proportion of the problem instances without a feasible solution at timeout over the 20 runs, with the time limit set to 1 hour.

Table~\ref{B-tab:results-time} shows that the IP model often times out without having found a feasible incumbent solution. This happens in all setups except the easiest, when $\lambda=16, W=2$. When $\lambda=18, W=4$, this occurs almost all the time for the IP model for objective \#1 and \#4, and also in a non-negligible proportion of instances for the CP model and the combined CP/IP approach. 

The results in Table~\ref{B-tab:results-quality} show that the CP model has the worst performance in almost all setups with regard to relative optimality gap after 1 hour. This can be expected: both \cite{Frimodig2019} and \cite{Pham2021} showed that CP is good at finding feasible solutions for the RT scheduling problem, but not as efficient at finding an optimal solution. For all objectives and setups, the CP model frequently times out without proving optimality within the time frame. However, columns \textbf{A-B} in Table~\ref{B-tab:results-quality} show that the relative optimality gap is often small.

For $\lambda=16$, the IP model often outperforms the Combined CP/IP model, suggesting the IP solver does not benefit from being warm started with a feasible CP solution in these cases. For more complicated instances ($\lambda=18$), the results are the opposite; the Combined CP/IP methodology gives shorter average CPU time and fewer instances that time out without having found a feasible solution than the pure IP model.

When altering the number of time windows, Table~\ref{B-tab:results-quality} shows that $W=2$ gives smaller relative optimality gaps and fewer instances without feasible solutions before timeout than for $W=4$. Moreover, Table~\ref{B-tab:results-time} shows the average solution times are also much shorter when $W=2$ than when $W=4$, likely due to the smaller problem dimensions. This difference in performance is largest for the IP model and smallest for the CG-IP model. 

Overall, the CG-IP approach has the smallest variance in solution times. Thus, the CG-IP model seems more robust to model size changes than the other models. This can also be seen when increasing the arrival rate from $\lambda=16$ to $\lambda=18$, as this increase has the smallest effect on the CG-IP model. Furthermore, the solution times are shortest for the CG-IP model, and the quality of solution is the best for this model.

\subsection{Objective function evaluation}
Evaluating the different objective functions, the results in Tables \ref{B-tab:results-time}-\ref{B-tab:results-quality} suggest that for the IP model, objective functions \#1 and \#4 are much harder to solve than objective functions \#2 and \#3. This can be seen both for the mean relative optimality gap, the proportion of instances with no feasible solution and the CPU times. In objective function \#1 and \#4, objective \ref{B-enum:6} is active, i.e., minimization of the number of switches to partially beam-matched machines. This seems to make the IP model much more complex, and the time for solving the root node relaxation alone increases from $200-300$ seconds for objective function \#2, to $700-800$ seconds for objective function \#4, although the problem dimensions are approximately the same. 

\begin{figure*}[b]
\center
\includegraphics[width=0.98\linewidth]{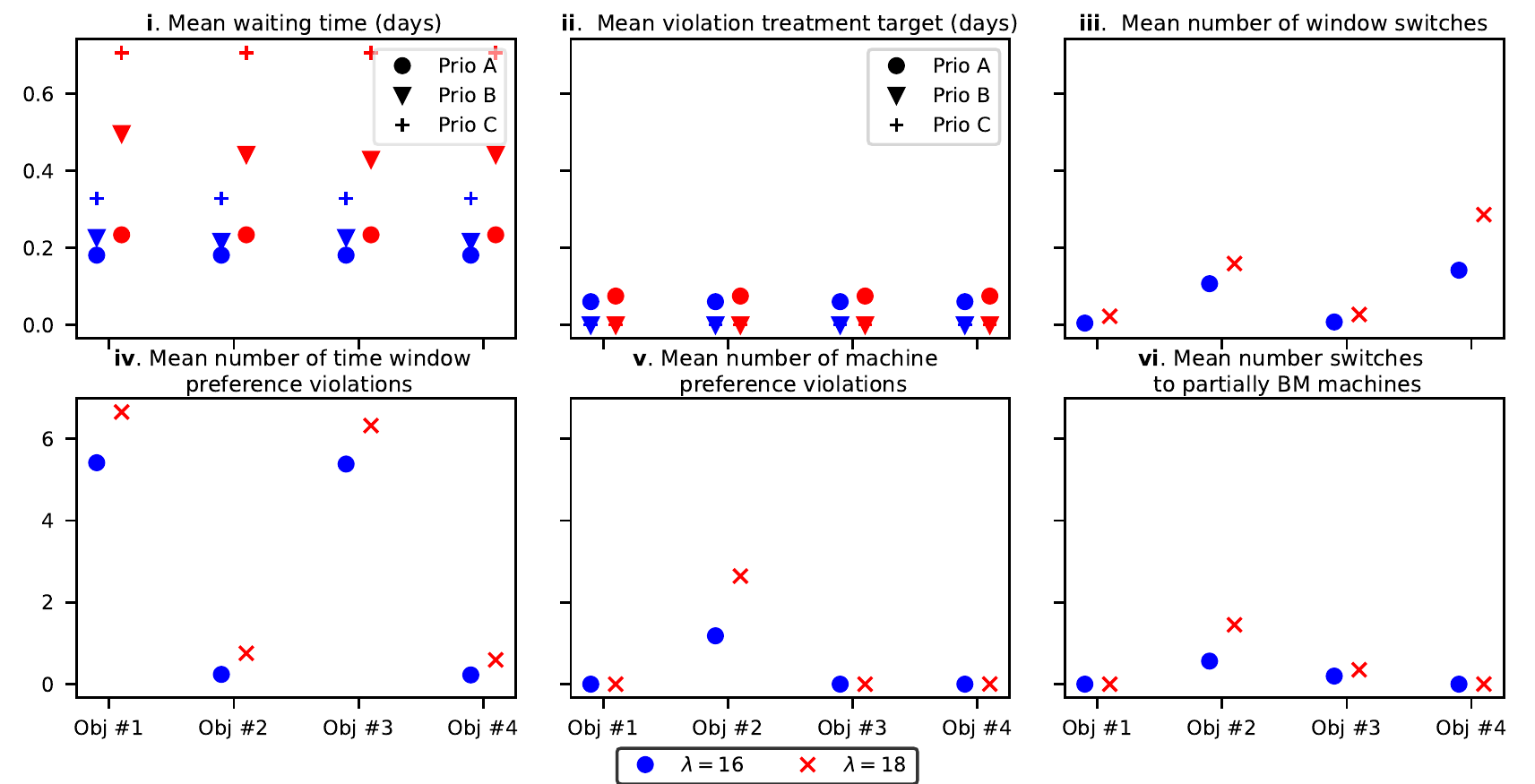}
\caption{Measurements from the CG-IP solutions relating to objective \ref{B-enum:i}-\ref{B-enum:6} in Table \ref{B-tab:obj_funs} are shown for the different objective functions when $W=4$} 
\label{B-fig:ObjPlot}
\end{figure*}

For the CP model, it is instead objective functions \#2 and \#4 that are more difficult than objective functions \#1 and \#3. In objective functions \#2 and \#4, objective \ref{B-enum:iv} is active, i.e., minimization of time window preferences. This objective makes it more difficult for the CP search heuristic to find the optimal solution, or even a feasible solution within the time limit. 

The CG-IP model is the only one that never times out without having found a feasible solution. This is expected, since the initial schedules generated by the heuristic in Algorithm~\ref{B-Alg:ColGenGreedy} are all feasible. The time results in Table~\ref{B-tab:results-time} show that the CG-IP model is also less sensitive to objective function changes than the other models. The results indicate that objective function \#1 may be less complicated since its solution times are shorter, but the difference to the other objective functions is smaller than the differences between objective functions for the other methods.

To evaluate the weights of the objective functions presented in Table~\ref{B-tab:obj_funs}, Figure~\ref{B-fig:ObjPlot} shows different costs from the CG-IP solutions for the different objective functions and different arrival rates. Each subplot \textbf{i}. to \textbf{iv}. represents the key indicator of each of the objectives \ref{B-enum:i} to \ref{B-enum:6}. From plot \textbf{i}. to the upper left, one can see that the mean waiting time does not change much between the different objective functions, which is the expected result since objective \ref{B-enum:i} is present in all objective functions. However, the waiting times increase as the arrival rate increases from $\lambda=16$ to $\lambda=18$ patients per day. The plot of the mean violations of treatment target times, \textbf{ii}., shows that this violation is always close to zero. This indicates that the addition of objective \ref{B-enum:ii} (to minimize the violation of treatment target dates) does not have a large, or any, effect when simultaneously minimizing the waiting times. 

To minimize the number of window switches, objective \ref{B-enum:iii}, is present in all objective functions. Subplot \textbf{iii}. in Figure~\ref{B-fig:ObjPlot} shows the mean number of window switches, and this value is very low for all objective functions. Objective \ref{B-enum:iv}, to minimize the violation of the time window preferences, is active in objective function \#2 and \#4. Subplot \textbf{iv}. shows that this is reflected in the results; this violation is much higher in objective function \#1 and \#3. To minimize the machine preference violations, objective \ref{B-enum:5}, is present in \#1, \#3 and \#4, which agrees with the results in subplot \textbf{v}. Finally, objective \ref{B-enum:6}, to minimize the number of switches to partially beam-matched machines, is present in objective function \#1 and \#4, and although subplot \textbf{vi}. shows that the mean value for the number of switches is low also for \#2 and \#3, it is lower for \#1 and \#4.

In total, this shows that the weights for the objective functions presented in Table~\ref{B-tab:obj_funs} are well reflected in the resulting schedules computed by the CG-IP model. It also shows that when the capacity is more limited due to a higher arrival rate, all the objectives are more difficult for the model to achieve.

\subsection{Sensitivity Analysis}

\begin{figure*}[b]
\center
\includegraphics[width=0.98\linewidth]{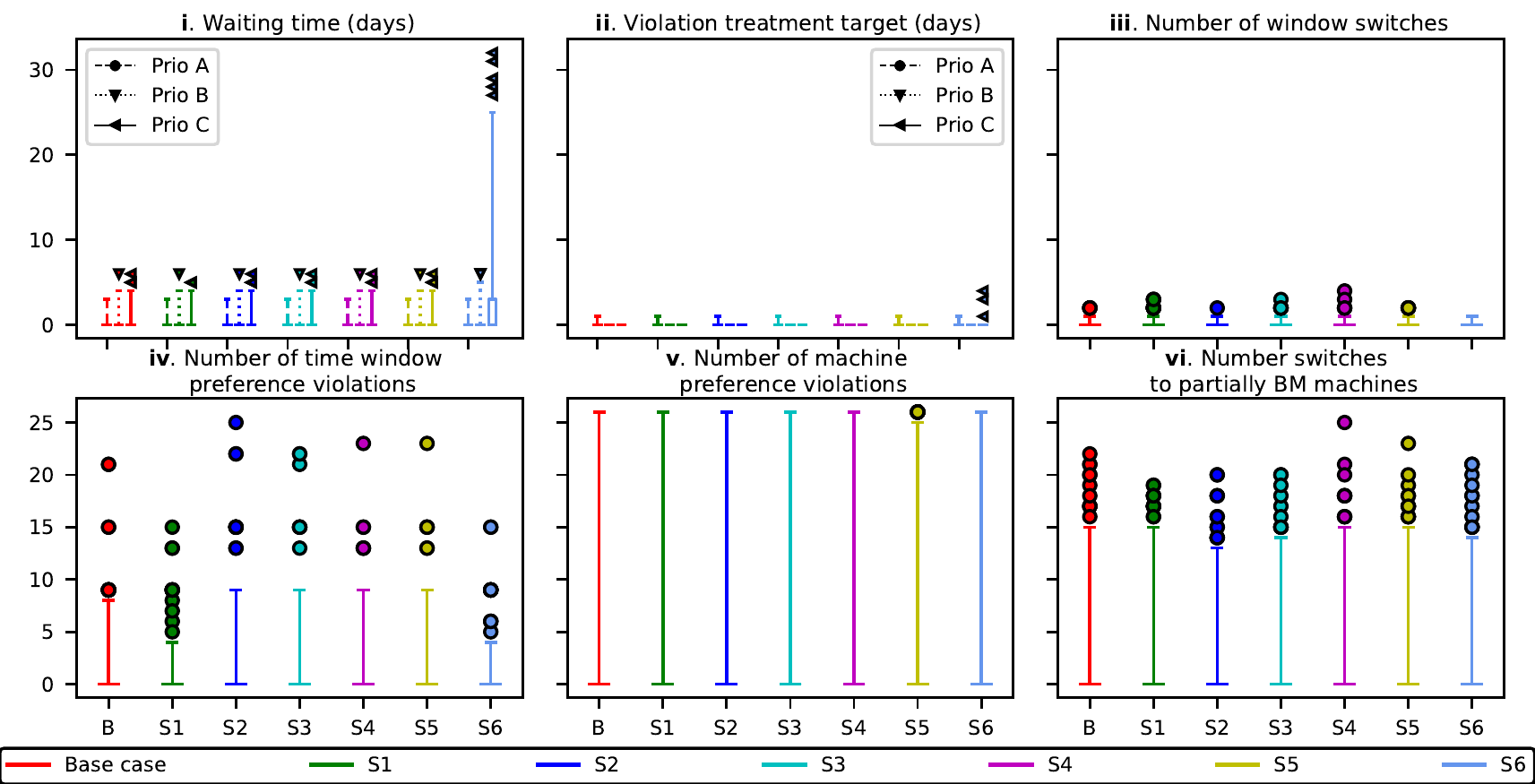}
\caption{Boxplots of the results in the different objectives in the sensitivity analysis using the CG-IP model. The weights of $\alpha$ for objective function \#2 are varied according to Table~\ref{B-tab:sensitivity}. The top 1\% are shown as outliers}
\label{B-fig:Sensitivity}
\end{figure*}

The parameters $\alpha_1,\dots,\alpha_4$ are included in the sensitivity analysis. Objectives \ref{B-enum:5} and \ref{B-enum:6} (relating to machine preferences and machine switches) are not relevant for clinics with a homogeneous machine setup. Therefore, the sensitivity analysis is focused on objective function \#2, for which $\alpha_5=\alpha_6=0$.


In Table~\ref{B-tab:sensitivity}, the base case used in the previous experiments and setups S1--S6 are presented. Because of the medical consequences, it is always more important to minimize waiting times for treatment start \ref{B-enum:i} (thereby also minimizing the violations of the treatment time targets \ref{B-enum:ii}) than to achieve a better patient experience  \cite{VanLent2013,Harden2022}, in this case by maximizing the time consistency in treatments \ref{B-enum:iii} and minimizing the violations of the patient wishes on treatment times  \ref{B-enum:iv}. This is reflected in the weights in the sensitivity analysis for all cases but one. In S6, it is instead prioritized to minimize the violation of the waiting time targets, and secondly to fulfill time consistency between appointments and to fulfill the patients time window preferences. Thirdly, the waiting times should be minimized. 

 \begin{table}
\caption{Sensitivity analysis for objective function \#2, where $\alpha_5=0, \alpha_6=0$}
\label{B-tab:sensitivity}
\begin{tabular}{ P{0.28\linewidth} P{0.62\linewidth} }
\hline\noalign{\smallskip}
Case & Weights \\
\noalign{\smallskip}\hline\noalign{\smallskip}
Base case  & $\alpha_1=50$, $\alpha_2=100$, $\alpha_3=1$, $\alpha_4=1$  \\
S1: Sensitivity 1 &$\alpha_1=10$, $\alpha_2=100$, $\alpha_3=1$, $\alpha_4=5$\\
S2: Sensitivity 2 & $\alpha_1=10$, $\alpha_2=100$, $\alpha_3=5$, $\alpha_4=1$ \\
S3: Sensitivity 3 & $\alpha_1=50$, $\alpha_2=50$, $\alpha_3=1$, $\alpha_4=1$\\
S4: Sensitivity 4 & $\alpha_1=100$, $\alpha_2=10$, $\alpha_3=1$, $\alpha_4=1$\\
S5: Sensitivity 5 & $\alpha_1=5$, $\alpha_2=10$, $\alpha_3=1$, $\alpha_4=2$\\
S6: Sensitivity 6 & $\alpha_1=1$, $\alpha_2=50$, $\alpha_3=5$, $\alpha_4=5$\\
\noalign{\smallskip}\hline
\end{tabular}
\end{table}

The results for the different parameter settings using the CG-IP model are shown in Figure~\ref{B-fig:Sensitivity}. Since objectives \ref{B-enum:5} and \ref{B-enum:6} are inactive in objective function \#2, it can be expected that the results do not differ very much between the parameter setups, which is confirmed by the results. 
Furthermore, the waiting times and the violations of the treatment time targets are almost identical between the different parameter setups if excluding S6. The mutual order of $\alpha_1$ and $\alpha_2$ does not seem to matter as long as both are greater than $\alpha_3$ and $\alpha_4$: in S4,  $\alpha_1$ has a higher weight than $\alpha_2$, which does not change the results in objectives \ref{B-enum:i} or \ref{B-enum:ii}. 
The results for objective  \ref{B-enum:iii} are similar for S1-S5, with differences only in the top 1\% shown as outlier points, except for S6. 
In the results for objective  \ref{B-enum:iv}, S1 and S6 have a lower number of time window preference violations than the other parameter setups. This is likely caused by the objective weight $\alpha_4$ being higher relative to $\alpha_3$ than in the other setups. Since the other results are very similar for S1 in the other metrics, this shows that for a clinic with a different prioritization between the objectives, it is possible to adjust the weights to achieve the required order. 

In S6, minimizing waiting time is no longer prioritized over minimizing time window switches and time window preference violations. This is probably not a relevant clinical scenario, but can give some insights in how the composite objective function works. The results in Figure~\ref{B-fig:Sensitivity} show that both the number of time window switches \ref{B-enum:iii}, and the number of time window preference violations \ref{B-enum:iv} are indeed lower for this setup, however, at the cost of some very long waiting times  \ref{B-enum:i}. 

The solution times for the different parameter setups are very similar to the base case times presented in Table~\ref{B-tab:results-time}. Overall, the sensitivity analysis shows that the composite objective function is not sensitive to the choice of the weights $\alpha_1-\alpha_6$, but their relative size order matters for the results.

\subsection{Conflicting Objectives}

\begin{figure}[h]
\center
\includegraphics[width=0.9\linewidth]{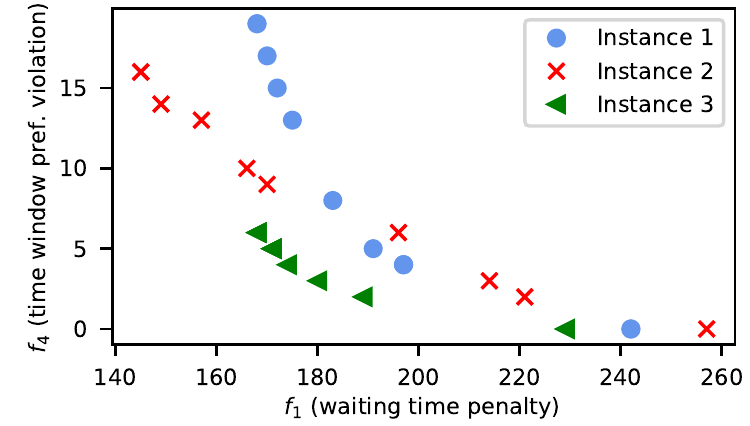}
\caption{Pareto optimal points for displaying the trade-off between \ref{B-enum:i} and \ref{B-enum:iv} for three instances where $\lambda=18$, $W=4$}
\label{B-Fig:pareto}
\end{figure}

The results of case S6 in the sensitivity analysis in Figure~\ref{B-fig:Sensitivity} indicate there could be a conflict between objectives \ref{B-enum:i} and \ref{B-enum:iv}, i.e., to minimize waiting time and to minimize the violations of the patients' time window preferences. 
Using the IP model and the weighted sum method for multi-objective optimization \cite{Zadeh1963,Marler2010}, this is analyzed for three randomly selected problem instances when $\lambda=18$, $W=4$. Figure~\ref{B-Fig:pareto} shows the pareto optimal points for the three instances when optimizing only the objectives \eqref{B-Eq:f1} and \eqref{B-Eq:f4} (both summed over $\pinP$). It shows that there is indeed a conflict between the objectives  \ref{B-enum:i} and \ref{B-enum:iv}; if only minimizing the waiting times, there are more violations of the time window preferences, and if minimizing the violations of time window preferences to optimality, the waiting time penalty will be much higher. 
Since there are severe medical consequences for having a longer waiting time, the trade-off between these two objectives is easily managed; the weight $\alpha_1$ for objective \ref{B-enum:i} should always be some magnitudes larger than $\alpha_4$ for objective \ref{B-enum:iv}.

\subsection{Clinical Implications}
The majority of the problem instances represent realistic scenarios since they are generated from clinical data. Objective function \#4 is most similar to what is used at Iridium Netwerk today, and $\lambda=16$ represent the average arrival rate at Iridium Netwerk. Therefore, this setup is used to analyze the clinical implications of using the CG-IP method for automatic scheduling. 

\begin{figure}[h]
\includegraphics[width=0.999\linewidth]{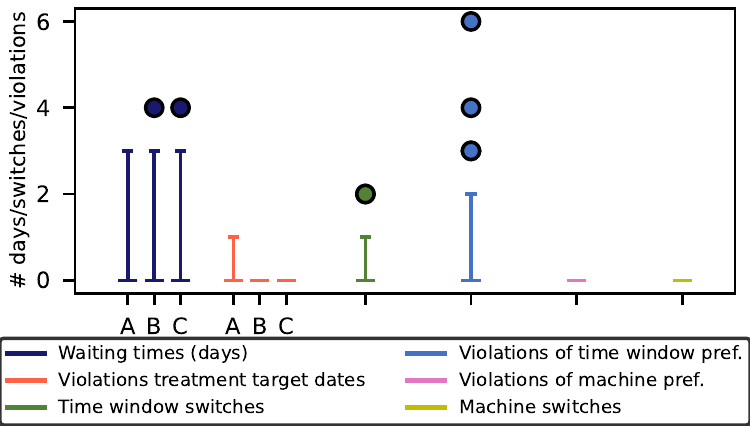}
\caption{Performance metrics for CG-IP for $\lambda=16$, $W=4$, with the top 1\% shown as outlier points}
\label{B-Fig:performance}
\end{figure}

Figure~\ref{B-Fig:performance} shows the results for the 20 problem instances where $\lambda=16$, $W=4$, for objective \#4, where a total of 327 patients have been scheduled to start treatment. Each of the six objectives described in Section~\ref{B-Sec: Problem Formulation} has its own boxplot, and for the waiting times and violations of waiting time targets, the results are further divided by the priority groups. 

The results demonstrate that the waiting times are always below one week. The exact clinical waiting times are not available, but the staff at Iridium Netwerk certify that they are often 2-4 weeks for priority B and C patients. These preliminary results show that there could potentially be large clinical benefit of using an the CG-IP model for automatic schedule generation. 

\section{Discussion}
\label{B-Sec:Discussion}
This section discusses the computational results, followed by a discussion of the potential for clinical implementation and directions for future work. 

\subsection{Model performance}
\label{B-Sec:Disc_performance}
The results show that the CG-IP model outperforms all other approaches in every aspect. It always finds feasible solutions and has the lowest mean optimality gap after one hour of run time. Table~\ref{B-tab:results-time} shows these solutions are almost always found long before the time limit; when there is a nonzero mean optimality gap, it is most often because the optimal solution has not been generated by the column generation procedure. This can happen since the CG algorithm (Algorithm~\ref{B-Alg:ColGen}) does not guarantee the optimal solution to be found. Table~\ref{B-tab:results-quality} shows the mean deviation from the optimal value is always below $1\%$, which means the solutions are of very good quality. Furthermore, the solution times of the CG-IP algorithm can possibly be decreased by solving $200-300$ independent subproblems in parallel. 

The CP model is the slowest of all models, and frequently times out without having found a feasible solution. For the smallest cases, when $\lambda = 16$ and $W=2$, the quality of the solution is very good although it often reaches the time limit. For $\lambda=16$ and $W=4$, it performs well for objective \#1 and \#3. The CP model could therefore be considered suitable for a clinical implementation if the clinic's workload is not too high, and especially if the clinic does not try to fulfill the patients' time window preferences. 

The IP model performs very well when $\lambda=16$ and $W=2$. Both when $\lambda = 16, W=4$ and when $\lambda=18, W=2$ the IP model also performs well for objective function \#2 and \#3. If a clinic does not need to support partially beam-matched machines, the IP model could therefore be suitable for clinical implementation. However, it is not suited for Iridium Netwerk, or other clinics where specific machine switches is an objective to be minimized. 

The results for the combined CP/IP methodology are better than the pure IP model, with fewer timeouts and better quality solutions. The disadvantage of developing and maintaining two separate models is however significant. Every time a constraint were to be altered or added, it would have to be done for both models, which also means trying out different formulations to optimize performance, and possibly change the CP search heuristic. Although the combined CP/IP methodology may not be maintainable in practice, it shows the potential for warm starting the IP model for difficult cases. A similar advantage could possibly be gained by warm starting the IP from the feasible schedule computed when calculating the time horizon (see Section~\ref{B-sec:time horizon}), which would require fewer computations and no extra model.

\subsection{Potential for clinical implementation}
\label{B-Sec:Disc_results}

The CG-IP model is considered most suitable for clinical implementation since it gives good quality solutions in a reasonable time, and is most robust to different setups. Robustness is crucial if the model should be generally applied to RT centers of different sizes. 

The clinical staff does not necessarily know anything about mathematical programming. Thus, for a clinical implementation a user interface is needed, where the computations can be performed in the background. For the CG-IP to work in a generalized clinical setting, a number of constraints and objectives should be implemented in the model, with the possibility to activate them in the used interface by a simple click. Furthermore, different trade-offs between the objectives should be available depending on the clinic's needs, which will correspond to alterations in the $\alpha$ values.


The models have been developed in collaboration with Iridium Netwerk, but the results of the algorithms should be generalizable. The majority of the constraints presented in Section~\ref{B-Sec: Problem Formulation} are applicable at other RT clinics as well, such as treatments on consecutive weekdays, specific allowed start days, specific machines suited for the treatments, different fraction durations, and machine capacity limits. When an objective is inactive in an objective function (by setting that particular $\alpha_i=0$), the variables relating to that objective are free to take any values and do not contribute to the solutions. 
It is possible that a clinic has some other medical or technical constraint that is not implemented in the models, but related literature (e.g. \cite{Conforti2010,Saure2012,Legrain2015,Vieira2021,Pham2021}) show that the majority of the collaborating RT centers have similar constraints. The largest difference is likely the scheduling strategy; some clinics require the patients to be scheduled immediately at arrival, which would not work with the presented models as they all require multiple patients in a batch to be scheduled.

The objective function evaluation and the sensitivity analysis also indicate that the algorithms are generalizable: if a clinic has a different prioritization order among the objectives, it is straightforward to change the weights to acquire the preferred order, and the as long as the weights differ in some orders of magnitude, the models are not sensitive to their exact values.

\subsection{Future work}
\label{B-Sec:future}

The models have so far been compared to each other. To further evaluate the CG-IP approach, the automatically generated schedules should be compared to manually constructed schedules from Iridium Netwerk. The schedules obtained from the models also need to be assessed for their practical feasibility by the clinical staff. 

There are various obstacles to implementing an automatic scheduling algorithm clinically. The first has been discussed above; a user interface is needed for a clinical implementation. Another is that the models do not support non-conventional treatments, such as multiple fractions per day, or treatments on non-consecutive days. At Iridium Netwerk, these are approximately 10\% of all treatments, so for an automated scheduling algorithm to work in practice this extension would be necessary. Furthermore, sometimes a linac is unavailable because of maintenance or software upgrades or holidays, which is not taken into account in the models. When modeling machine unavailability, constraints on the minimum number of fractions per week are important for the treatment outcome. Finally, the patient protocols and the patient preferences regarding treatment time of the day would have to be registered instead of communicated verbally. 

The method of using placeholder patients to account for expected future patients should be evaluated. Perhaps it is adequate to reserve time on the linacs each day for future patients instead, which would decrease the solution times. If not, a potential improvement of the stochastic aspect of the models is to use scenario-based probabilities instead of expected values. Another method would be to use a data-driven approach to predict future machine utilization.

Since there are multiple objectives to be optimized, multi-objective optimization could be considered. It is possible that a multi-objective approach could be relevant at some clinics, but not at Iridium Netwerk; one of the main objectives with automatic schedule creation is to minimize the manual work, and there is no desire from the clinical staff to be able to navigate trade-offs in several generated schedules. Furthermore, the computational time to generate the schedule cannot be too long for it to work in clinical practice, which makes multi-objective optimization unfit for the task.

A future extension of the automatic scheduling approach would be to be able to optimize schedules for the whole appointment series of each patient, including not only the linac scheduling, but also meetings with physicians and RTTs, the duration of the treatment planning, CT scans, and more. It would be interesting to study the trade-off between treatment plan quality and waiting time - when is it beneficial for the patient to have a longer waiting time in order to be treated on a more advanced machine, and when is a less advanced machine with shorter waiting time to prefer? Perhaps the toxicity effect of assigning patients to non-preferred machines could be incorporated in the scheduling models. To optimize the use of resources and the treatment plans simultaneously could potentially have a great impact on both the clinics and the individual patients.

\section{Conclusions}
\label{B-Sec: Conclusions}

As the incidences of cancer increase, the demand for RT grows. To better use resources in RT, algorithms can be used to automatically create patient schedules, a task that today is done manually in almost all clinics. The main contribution of this paper is to serve as a decision support tool when implementing a scheduling algorithm in practice. We present a comprehensive study of a wide range of optimization methods that can be used to model the RT scheduling problem. The output of the models is the assignment of all fractions of the patients to both linacs and specific time windows, while including all the constraints necessary for the scheduling to work in practice. The models developed include an IP model, a CG-IP model, and a CP model, as well as a method combining the CP and IP models.

The models are tested using historical data from Iridium Netwerk in Antwerp, Belgium. Different cancer centers may have different intentions when creating the RT schedules, and in order to study the suitability of the different models for various cancer centers, each model is solved using multiple different objective functions. This is to evaluate if some particular optimization method is better suited to solve a certain objective. 

The results demonstrate that the CG-IP model is the most robust, and that the mean optimality gap of the method is well below $1\%$ for all the different setups and objective functions after one hour of computation time. The CP and IP models could have potential for clinical implementation depending on the size of the clinic, and more importantly, depending on their objective of scheduling.

The proposed methodology provides a tool for automated scheduling of RT treatments on linacs, and can be generally applied to RT centers. This would allow the RT staff to save time, and at the same time create optimized patient schedules that take medical and technical constraints into account. Designing more efficient schedules could potentially save lives by shortening waiting times and improving patient outcomes. 

\section*{Conflicts of Interest}
\vspace{-7pt}
On behalf of all authors, the corresponding author states there is no conflict of interest.

\label{Bibliography}
\bibliographystyle{unsrt} 
\bibliography{mybibliography}

\end{document}